\documentclass[12pt]{article}
\usepackage{amsfonts, amssymb, amsmath}

\usepackage{enumerate} 

\usepackage{bbm} 

\usepackage{scalerel}

\usepackage{color}
\definecolor{lightblue}{rgb}{0,0.2,0.5}

\usepackage[colorlinks=true, urlcolor=blue,linkcolor=blue, citecolor=lightblue]{hyperref}

\usepackage[round]{natbib}

\bibpunct{\textcolor{lightblue}{(}}{\textcolor{lightblue}{)}}{,}{a}{}{;}

\usepackage{graphics,graphicx,amsmath}

\newcommand{\R}{\mathbb{R}}

\newcommand{\E}{\mathbb{E}}

\newcommand{\N}{\mathbb{N}}

\let\Horig\H

\newtheorem{prop}{Proposition}[section]
\newtheorem{lemma}[prop]{Lemma}
\newtheorem{definition}[prop]{Definition}
\newtheorem{corollary}[prop]{Corollary}
\newtheorem{theorem}[prop]{Theorem}

\def\P{{\mathord{\mathbb P}}}

\def\({\left(}
\def\){\right)}

\def\[{\left[}
\def\]{\right]}
\def\real{{\mathord{\mathbb R}}}
\def\inte{{\mathord{\mathbb N}}}
\def\Dom{\mathrm{Dom}}
\def\Var{\mathrm{Var}}

\newenvironment{Proof}{\removelastskip\par\medskip
\noindent{\em Proof.} \rm}{\penalty-20\null\hfill$\square$\par\medbreak}

\allowdisplaybreaks

\numberwithin{equation}{section}

\begin{document}

\title{
\huge
Berry-Esseen bounds for functionals of independent random variables 
} 

\author{Nicolas Privault\thanks{Division of Mathematical Sciences,
Nanyang Technological University, SPMS-MAS-05-43, 21 Nanyang Link
Singapore 637371. e-mail: {\tt nprivault@ntu.edu.sg}. 
}
\and 
Grzegorz Serafin\thanks{Faculty of Pure and Applied Mathematics, Wroc{\l}aw University of Science and Technology, Ul. Wybrze\.ze Wyspia\'nskiego 27, Wroc{\l}aw, Poland.
e-mail: {\tt grzegorz.serafin@pwr.edu.pl}.}}

\maketitle 

\vspace{-0.8cm} 
\begin{abstract} 
 We derive Berry-Esseen approximation bounds
 for general functionals of independent random variables, based on
 chaos expansions methods. 
  Our results apply to $U$-statistics satisfying the weak assumption of
decomposability in the Hoeffding sense, and yield Kolmogorov distance
bounds instead of the Wasserstein bounds previously derived in the
special case of degenerate $U$-statistics.
  Linear and quadratic functionals 
 of arbitrary sequences of independent  random variables are included as particular cases, 
 with new fourth moment bounds, and 
 applications are given to Hoeffding decompositions,
 weighted $U$-statistics, quadratic forms, and random subgraph weighing.
 In the case of quadratic forms, our results
 recover and improve the bounds available in the literature,
  and apply to matrices with non-empty diagonals.
\end{abstract}

\baselineskip0.7cm

\section{Introduction}
Significant progress in probability approximation has been
achieved in recent years by combining the Chen-Stein method 
with the Malliavin calculus.
 See for example \cite{nourdinpeccati}, \cite{utzet2}, \cite{thale},
for the derivation of distance bounds on the Wiener and Poisson spaces, 
and also \cite{nourdin3} and \cite{reichenbachs} 
 in the case of Rademacher sequences. 
 Those results rely on covariance
 representations based on the inverse 
 of the Ornstein-Uhlenbeck operator $L$ acting on multiple
 Wiener-Poisson stochastic integrals. 
 While the inverse operator $L^{-1}$ is well adapted to
 certain random functionals such as multiple stochastic integrals,
 it can prove 
 more difficult to use in applications to other, more specific functionals.
 Other covariance representations based on
 the Clark-Ocone representation formula and not relying on $L^{-1}$
 have been used in \cite{privaulttorrisi3} 
 on the Wiener and Poisson spaces, and in \cite{privaulttorrisi4}
 for Rademacher sequences. 

  \medskip

 In \cite{lastpeccatipenrose}, second order Poincar\'e inequalities
 in the
 Kolmogorov and Wasserstein distances 
 have been obtained for functionals of a Poisson
 point process by using the iterated Malliavin gradient instead of $L^{-1}$.
This approach relies on probabilistic representations for
the inverse operator $L^{-1}$ using
Mehler's formula on the Poisson space,
see e.g. Lemma~6.8.1 in \cite{privaultbk2}.
Second order Poincar\'e inequalities for functionals of
Rademacher sequences have also been obtained in \cite{krokowski},
with application to renormalized triangle counting using the Kolmogorov 
distance in the Erd{\H o}s-R\'enyi random graph,
see also \cite{PS2} and references therein for the treatment of
arbitrary subgraph counting. 

\medskip
 
In \cite{PS}, a general framework for the
derivation of Wasserstein distance bounds
for functionals of independent random sequences
has been developed in the integration by parts setting
of \cite{prebub}, using an analog of the
operator $L^{-1}$ on discrete chaos expansions based on discrete multiple
stochastic integrals.
Bounds in total variance distance have also been obtained
therein using Clark-Ocone covariance representation formulas
under stronger smoothness conditions. 
Related results have been obtained in \cite{halconruy}, see Theorem~5.9
therein for a normal Stein approximation bound for functionals of
independent random variables, see also \cite{tdnguyen}, and \cite{sambale}
for concentration inequalities. 
Applications to normal approximation
in the Wasserstein distance have been obtained in \cite{PS3} 
for the weights of subgraphs in the Erd{\H o}s-R\'enyi random graph. 
 
 \medskip

 Our first goal in this paper is to extend existing
 Stein normal approximation bounds proved in the Kolmogorov distance
 for Rademacher sequences, see e.g.
 \cite{krokowski}, \cite{dobler}, 
 to general sequences of independent random
 variables. 
 This is achieved in the general framework of \cite{PS}, 
 by replacing the Wasserstein distance with the Kolmogorov distance  
 for which obtaining rates is known to be more difficult and
 requires new ideas.
In Theorem~\ref{thm:main} we derive a general Berry-Esseen bound 
which is then specialized to
sums of multiple stochastic integrals in
Proposition~\ref{prop:dK<Sum}
and then to multiple stochastic integrals in
Proposition~\ref{prop:dKI}.
 Note that multiple stochastic integrals of order $d$ coincide with  degenerate (generalized) $U$-statistics of order $d$, and can then
  be used to represent Hoeffding decompositions
  as a chaos summations, see the examples given below.

\medskip

{Our second goal is to show that the obtained bounds remain sharp
 despite the very general framework of the paper,
  as demonstrated in the following examples}.
 Consider a sequence $(X_1, \ldots , X_n)$ of
(not necessarily identically distributed) independent random variables,
 and the $d$-homogeneous random multilinear forms $W_{n,d}$
 written in the Hoeffding form as 
$$
W_{n,d} = \sum_{J \subset \{1,\ldots , n \}, \ |J|=d} W_J, 
$$
where, for each $J \subset \{1,\ldots , n \}$, 
$W_J$ is a random variable with variance $\sigma^2_J$, 
measurable with respect to the $\sigma$-algebra 
$\mathcal F_J :=\sigma\(X_j:j\in J\)$,
and such that $\E\[W_J \mid \mathcal F_K\]=0$, 
$J\not\subseteq K\subset [n]$. 
 In \cite{dejong1990}, 
a central limit theorem has been proved 
for the sequence $(W_{n,d})_{n\geq 1}$ under the conditions
$$
\lim_{n\to \infty} \max_{1\leq i \leq n} \sum_{J\ni i} \sigma_J^2 = 0
\quad
\mbox{and}
\quad 
\lim_{n\to \infty} \E [ W_{n,d}^4 ] = 3,
$$
 generalizing earlier results by \cite{dejong1987}
 for quadratic random functionals.
 The results of \cite{dejong1987,dejong1990} 
 have been refined 
 by the derivation of bounds in the Wasserstein distance in
 Theorem~1.3 in \cite{doblerpeccati} in the case 
 of degenerate $U$-statistics, 
   for which $|J|$ is constrained to a fixed value $|J|=d$ for
   some $d\in \{1,\ldots , n\}$ in the sum
   \eqref{eq:Hd}. 

\medskip 

Applications of Proposition~\ref{prop:dK<Sum} 
 are given to Kolmogorov distance bounds 
 in Theorem~\ref{thm:S}
 for general $U$-statistics, and in Theorems~\ref{thm:degenerateU} and \ref{thm:Ustat} for degenerate $U$-statistics.
 This extends 
 the bounds of \cite{doblerpeccati} by using 
 the Kolmogorov distance instead of the Wasserstein
   distance, and by applying to Hoeffding decompositions
   in full generality and not only to degenerate $U$-statistics. 
 This also extends the bounds in the Kolmogorov distance 
 derived in \cite{dobler} for $U$-statistics
 in the particular case of Rademacher chaoses,
 where $(X_1,\ldots , X_n)$ is a sequence of independent 
 Bernoulli random variables. 

   \medskip

 More specifically, given an i.i.d.
sequence $(X_k)_{k\geq 1}$ of centered random variables 
with unit variance,
 and the sum 
$$
Z_n : = \frac{1}{\sqrt{n}} \sum_{k=1}^n X_k, \qquad n \geq 1, 
$$
 convergence bounds to the standard normal distribution 
$\mathcal N$ of the form 
\begin{equation}
\nonumber 
 d_W( Z_n ,\mathcal N) \leq \frac{E\[|X_1|^3\]}{\sqrt{n}} 
\end{equation} 
 have been obtained in e.g. Theorem~1.1 in \cite{goldstein} 
 in the Wasserstein distance 
$$
 d_W (X,\mathcal N ):
 =\sup_{h\in\mathrm{Lip}(1)}|\E [h(X)]-\E [h(\mathcal N )]|. 
$$
 See also Corollary~2.11 of \cite{dobler2} 
 for related bounds in the Kolmogorov distance 
$$
 d_K ( X,{\mathcal N}) : = \sup_{x\in \real}
 | P(X\leq x ) - P(\mathcal{N} \leq x ) |,
$$
 including the case of random sums. 
 In the case of quadratic functionals of the form
\begin{equation}
   \label{qn} 
 Q_n : = \sum_{1 \leq k , l \leq n} a_{kl} X_k X_l,
\end{equation} 
 where $A= ( a_{ij} )_{1\leq i,j\leq n}$ is a symmetric matrix, 
 the bound 
 \begin{equation}
   \label{gt1} 
 d_K ( Q_n , \mathcal N)\leq C \(E\[X_1^3\]\)^2 | \lambda_1 |, 
\end{equation} 
 where $\lambda_1$ denotes the largest absolute eigenvalue of $A$ 
 and $C>0$ is an absolute constant,
 has been obtained in \cite{goetze} when the diagonal of $A$ vanishes, 
 see e.g. Theorem~1 therein, 
 and also 
 Theorem~3.1 of \cite{shao2}. 

 \medskip

 In this 
 vanishing diagonals setting, Theorem~\ref{thm:Ustat} is applied to derive 
 Corollary~\ref{chkskl} which recovers Theorem~3.1 in \cite{shao2},   
 and improves on the above bound \eqref{gt1} of Theorem~1 in \cite{goetze}.
 In addition, Corollary~\ref{chkskl} extends 
 the Kolmogorov bounds of 
 Theorem~1.1 in \cite{dobler}, restricted to the quadratic case, 
 from Rademacher sequences to general sequences of random variables by 
 using fourth moment differences as in e.g. Theorem~1.3 of \cite{doblerpeccati}.  

 \medskip

 In case the diagonal of $A= ( a_{ij} )_{1\leq i,j\leq n}$ 
 may not vanish, the bound 
\begin{align}\label{eq:GT0}
 d_K\(\frac{Q_n}{\sigma_n},\mathcal{N}\)\leq C(\gamma)\frac{\(\E\[|X|^3\]\)^2+\gamma \E[X^6]}{\sqrt{\sum_{1\leq i,j\leq n}a_{ij}^2}}|\lambda_1|,
\end{align}
has been obtained in Theorem~1.1 of
\cite{goetze-tikhomirov2002} for some $\gamma >0$ depending on $A$.  
See also Proposition~3.1 in \cite{chatterjee} for a result in the Wasserstein
distance using Rademacher sequences, and Theorem~2.2 in \cite{chatterjee2} 
for related normal approximation bounds in total variation distance for a smooth
 function of finite-dimensional random vectors
 via second order Poincar\'e inequalities. 

 \medskip 
 
 In comparison with Theorem~1.1 of \cite{goetze-tikhomirov2002}, 
 the bound \eqref{r2} in Theorem~\ref{qf} gives better rates under weaker assumptions
 according to the inequality \eqref{eq:Tr<l}. 
 Theorem~\ref{qf} also provides an additional bound  
 \eqref{r1} which is valid for any i.i.d. sequence $(X_n)_{n\geq 1}$ 
 and holds in the Kolmogorov distance, instead
 of the Wasserstein distance used in \cite{doblerpeccati}.
 This bound is related to the so-called fourth moment phenomenon
 (\cite{nualart2005central}), which has been the object of intense research work,
 see e.g. \cite{Nourdin-Peccati-book} and references therein. 

\medskip
 
We proceed as follows.
In Section~\ref{s2} we recall the framework of 
\cite{prebub} for the treatment of functionals of
independent random sequences, including the
construction of discrete multiple stochastic integrals
and the associated finite difference gradient operator and 
integration by parts formula, which are used to derive a fourth moment bound
in Section~\ref{s2.1}. 
Section~\ref{s3} contains our main result Theorem~\ref{thm:main}
which states a general Berry-Esseen bounds
 for general functionals of independent random sequences.
In Section~\ref{s4} those results are applied to
the derivation of Kolmogorov bounds for discrete multiple
integrals and for sums of discrete multiple integrals. 
Applications to Hoeffding decompositions,
weighted $U$-statistics  and
subgraph counting in the 
 \cite{ER}
 random graph are given in Section~\ref{s5}.
 Section \ref{section:qf} focuses on quadratic forms.
 
\section{Preliminaries and notation} 
\label{s2}
 Consider an i.i.d. sequence $(U_k)_{k\geq 1}$ of uniformly distributed random
variables on the interval $[-1,1]$, 
on a probability space $(\Omega , {\cal F}, \P)$. 
Given $F(U_1(\omega),U_2(\omega), \ldots )$
a functional 
of the sequence $\(U_1(\omega),U_2(\omega), \ldots \)$,
we consider the shifted sequence $\Phi_t (\omega )$ defined as 
$$\Phi_t(\omega) : =\(U_1(\omega),\ldots ,U_{\lfloor t/2\rfloor}(\omega),t-1-2\lfloor t/2\rfloor,U_{\lfloor t/2\rfloor+2}(\omega),\ldots \),
\qquad t\in \real_+,
$$
and define the finite difference gradient operator $\nabla$
 on random functionals as 
\begin{equation}\label{eq:nablaaspsi}
  \nabla_tF := F \circ \Phi_t - \frac{1}{2}
  \int_{2\lfloor t/2\rfloor}^{2\lfloor t/2\rfloor+2}F\circ \Phi_s ds,
  \qquad t \in \real_+, 
\end{equation}
provided that
$( F\circ \Phi_s)_{s\in \real_+}$ is integrable on $\real_+$,
$\P$-a.s., 
see Definition~5 and Proposition~10 in \cite{prebub}.

For any $X \in L^1 ( \Omega )$ and $k \in \inte$
we also note the identity
\begin{equation}
  \label{ex} 
\E [ X ] = \frac{1}{2} \E \left[ \int_{2k}^{2k+2}
  X \circ \Phi_u du \right]. 
\end{equation} 
 Consider now the adjoint $\nabla^*$ of $\nabla$,
 defined by the duality relation 
\begin{equation}\label{eq:duality}
\E\[\langle\nabla X,u\rangle_{L^2(\R_+,dx/2)}\]=\E[X\nabla^*(u)], 
\end{equation}
 which shows that $\nabla^*$ and $\nabla$ are closable with domains
 $\Dom ( \nabla^* ) \subset L^2 (\Omega )$ and 
 $$
\Dom ( \nabla ) =
\big\{ X\in L^2(\Omega ) : E [ \Vert \nabla X \Vert^2_{L^2(\real_+)} ]
< \infty \big\}
\subset L^2 (\Omega \times \real_+ ),
$$
 see Proposition~8 in \cite{prebub}. 
 The operators $(\nabla, \nabla^* )$ are linked by the Skorohod isometry 
\begin{equation*}
  \E[\nabla^*u\nabla^*v]=\E\[\int_0^\infty u_tv_tdt\]+\E\[\int_0^\infty
  \int_0^\infty
  \nabla_su_t\nabla_tv_s\,ds\,dt\], 
\end{equation*}
see Proposition~9 in \cite{prebub},
which yields the Poincar\'e inequality
\begin{equation} 
  \label{eq:Enabla*^2}
  \E\big[ |\nabla^*u|^2\big] 
  \leq\E\[\int_0^\infty |u_t|^2dt\]+\E\[\int_0^\infty
  \int_0^\infty
  |\nabla_su_t|^2ds\,dt\]. 
\end{equation} 

\begin{definition}
  Given $f_n$ in the space $\widehat{L}^2(\R_+^n)$ of 
square integrable symmetric functions on $\R_+^n$ that
vanish outside of 
$$
\Delta_n : = \bigcup_{
  k_i \not= k_j \geq 1 \atop
  1\leq i \not= j \leq n
}
      [2k_1-2,2k_1]\times \cdots \times [2k_n-2,2k_n], 
$$
      we define the multiple stochastic integral
$$ 
 I_n(f_n)
  = 
  n!\int_0^\infty \int_0^{t_n}\cdots \int_0^{t_2}f_n(t_1,\ldots ,t_n)d(Y_{t_1}-t_1/2)\cdots d(Y_{t_n}-t_n/2), 
  $$
 with respect to the jump process 
$\displaystyle Y_t : =\sum_{k=1}^\infty\mathbf1_{[2k-1+U_k,\infty)}(t)$,
$t\in \real_+$, which satisfies 
\begin{align}  
& I_n(f_n)
   := \sum_{r=0}^n
    \left( - \frac{1}{2}\right)^{n-r} {{n}\choose{r}} 
   \\
  \nonumber
  & \qquad 
  \times \sum_{ k_1\neq\cdots \neq k_r \geq 1}
  \int_0^\infty\cdots \int_0^\infty f_n (2k_1-1+U_{k_1}
  ,\ldots ,2k_r-1+U_{k_r} ,y_1,\ldots ,y_{n-r})dy_1\cdots dy_{n-r}.  
\end{align} 
\end{definition}
 The multiple stochastic integral  
 $I_n(f_n)$ satisfies the bound 
 \begin{equation}
   \nonumber 
  \E \[( I_n(f_n))^2\]\leq n!\left\|f_n\right\|^2_{L^2(\R^n_+,dx/2)},
  \qquad n\geq 1, 
\end{equation} 
which allows us to extend
the definition of $I_n(f_n)$ to all
$f_n \in \widehat{L}^2(\R_+^n)$, see Propositions~4 and 6 in \cite{prebub}.
 Under the additional condition 
\begin{equation}
\label{ass:int0}
\int_{2k-2}^{2k}f_n (t,*)dt=0, \qquad k \geq 1, 
\end{equation} 
 i.e. $f_n$ is canonical in the sense of \cite{surgailisclt},
 the multiple stochastic integral $I_n(f_n)$
 can be written as the $U$-statistics of order $n$ 
\begin{align}
\nonumber 
 I_n(f_n)&=\sum_{ k_1\neq\cdots \neq k_n \geq 1 }f_n (2k_1-1+U_1,\ldots ,2k_n-1+U_n),
\end{align}
 with the isometry and orthogonality relation
\begin{equation}
\label{eq:EI^2}
\E\[ I_n(f_n) I_m(f_m) \]=
           {\bf 1}_{\{ n = m \} }
           n!\langle f_n , f_m \rangle_{L^2(\R_+,dx/2)^{\circ n}}, 
           \quad f_n \in \widehat{L}^2(\R_+^n), \quad
           f_m \in \widehat{L}^2(\R_+^m), 
\end{equation}
see \cite{prebub}, page~589,
which shows that 
the sequence $(I_n(f_n))_{n\geq 1}$ forms a 
family of mutually orthogonal centered random variables.
Finally, every $X\in L^2(\Omega)$ admits the chaos decomposition 
\begin{equation}\label{eq:decomp}
X=E[X] + \sum_{n=1}^\infty I_n(f_n),
\end{equation}
for some sequence of functions $f_n$ in $\widehat{L}^2(\R_+^n)$,
$n\geq 1$, cf. Proposition~7 in \cite{prebub}.
Moreover, under the condition \eqref{ass:int0}
the sequence $(f_n)_{n\geq 1}$ is unique
in $\widehat{L}^2 (\R_+^n)$ due to the isometry
relation~\eqref{eq:EI^2}, and in this case we have 
\begin{align}\label{eq:chaosnorm}
\E [X^2 ]=(\E[F])^2+\sum_{n=1}^\infty n!\|f_n\|^2_{L^2(\R_+^n,(dx/2)^{\otimes n})}.
\end{align}
 Under the condition \eqref{ass:int0} we also have the relations 
\begin{equation}
\label{alh} 
 \nabla^* \left( I_n(g_{n+1}) \right) : = I_{n+1}({\bf 1}_{\Delta_{n+1}}
 \tilde{g}_{n+1})
\quad
 \mbox{and}
 \quad 
  \nabla_t I_n(f_n) = n I_{n-1}\(f_n(t,*)\), \quad
 t\in \real_+,
\end{equation} 
 see Proposition~10 in \cite{prebub},
 where $\tilde{g}_{n+1}$ is the symmetrization of
$g_{n+1} \in \widehat{L}^2(\R_+^n)\otimes L^2(\R_+)$
in $n+1$ variables. 
The operator $L$ defined on linear combinations of multiple stochastic
integrals as
\begin{equation*}
  L I_n(f_n) : =-\nabla^*\nabla_t I_n(f_n) =-n I_n(f_n),
  \qquad
  f_n\in \widehat{L}^2(\R_+^n), 
\end{equation*}
is called the Ornstein-Uhlenbeck operator. By \eqref{eq:decomp}
the operator is invertible for centered
$X\in L^2(\Omega)$, and its inverse operator $L^{-1}$ is given by
\begin{equation}
  \label{L-1}
  L^{-1} I_n(f_n)  =-\frac{1}{n} I_n(f_n),\qquad n\geq1.
  \end{equation} 
  In fact, we can easily derive the form of any real power of $-L$, i.e. it holds
$$(-L)^{\alpha} I_n(f_n)  =n^\alpha I_n(f_n),\qquad n\geq1,\ \ \alpha\in\R.$$

 We also recall that, by Proposition~5.3 in \cite{PS3},
 for every $f_n\in \widehat{L}^2(\R_+^n)$
 there exists $\bar{f}_n\in \widehat{L}^2(\R_+^n)$
 given by
 \begin{align}
   \nonumber 
  \bar{f}_n(t_1,\ldots ,t_n)
= 
\Psi_{t_1} \cdots \Psi_{t_n} f_n(t_1,\ldots ,t_i),
\end{align}
  satisfying \eqref{ass:int0}
  and such that $I_n(f_n)=I_n (\bar{f}_n )$,
  where
\begin{align*}
  \Psi_{t_i}f(t_1,\ldots , t_n):=f(t_1,\ldots , t_n)-
  \frac{1}{2}
  \int_{2\lfloor t_i/2\rfloor}^{2\lfloor t_i/2\rfloor+2}f(t_1,\ldots ,t_{i-1},s,t_{i+1},\ldots ,t_n)ds,
\end{align*}
$i=1,\ldots , n$,
$t_1,\ldots ,t_n \in \real_+$.
We end this section with the following multiplication formula
for multiple stochastic integrals,
see Proposition~5.1 in \cite{PS}.
Letting $n\wedge m := \min (n,m)$, for $0\leq l\leq k\leq n\wedge m$
we define the contraction 
$f_n\star_{k}^lg_m$ of
$f_n \in \widehat{L}^2(\R_+^n)$ and
$g_m \in \widehat{L}^2(\R_+^m)$ as 
\begin{align}
  \nonumber 
 & f_n\star_{k}^lg_m(y_1,\ldots , y_{n-l},z_1,\ldots , z_{m-k}) 
  \\\nonumber
   & := \frac{1}{2^l} 
\int_{\R^l_+} f_n(x_1,\ldots , x_l ,y_1,\ldots , y_{n-l} ) g_m(x_1,\ldots , x_l ,y_1,\ldots , y_{k-l},z_1,\ldots , z_{m-k})dx_1\cdots dx_l,
\end{align}
           and we let
           $f_n \hskip0.1cm \widetilde{\star}_{k}^l g_m$
           denote the symmetrization
           \begin{eqnarray*}
             \lefteqn{ 
               \! \! \! \! \! \! \! \! \! 
               f_n \hskip0.1cm \widetilde{\star}_{k}^l \hskip0.1cm g_m(x_1,\ldots , x_{n+m-k-l})
             }
             \\
             & := &
             \frac{\mathbbm{1}_{\Delta_{m+n-k-l}}(x_1,\ldots , x_{n+m-k-l})
              }{(m+n-k-l)!}\sum_{\sigma\in \Sigma_{m+n-k-l}}f_n\star_{k}^lg_m(x_{\sigma(1)},\ldots ,x_{\sigma(m+n-k-l)}).
\end{eqnarray*} 
Then, for $f_n \in \widehat{L}^2(\R_+^n)$
and $g_m \in \widehat{L}^2(\R_+^m)$
satisfying \eqref{ass:int0}, the following multiplication formula holds: 
\begin{equation}
 \label{eq:mf}
 I_n(f_n) I_m(g_m)=\sum_{k=0}^{m\wedge n}k!{{m}\choose{k}}{{n}\choose{k}}\sum_{i=0}^k{{k}\choose{i}} I_{m+n-k-i}\big(
 f_n\hskip0.1cm \widetilde{\star}_{k}^{i} g_m\big),
\end{equation}
whenever $f_n\star_{k}^{i}g_m\in L^2(\R_+^{m+n-k-i})$
for every $0\leq i \leq k\leq m\wedge n$. 

\section{Fourth moment bound}
\label{s2.1}
The next covariance relation can be obtained as in Proposition~2.1 in \cite{hp}.
\begin{prop}
  \label{sdajklsad}
  Let $\alpha\in \real$
  and $X,Y\in L^2(\Omega )$ such that
  $L^{\alpha -1} X \in \Dom ( \nabla )$
  and
  $L^{- \alpha } Y \in \Dom ( \nabla )$.
  Then we have the covariance relation 
\begin{equation} 
  \label{eq:intnablaX}
  {\rm Cov}\(X,Y\)
= 
\frac12 \E\[\int_0^\infty (\nabla_t(-L)^{\alpha -1}X)(\nabla_t(-L)^{-\alpha} Y)\frac{dt}{2}\]. 
\end{equation} 
\end{prop}
\begin{Proof}
 We have 
\begin{eqnarray} 
\nonumber
  {\rm Cov}\(X,Y\)
  & = &
    \E [ ( X - \E [ X ] ) ( Y - \E [ Y ] ) ] 
    \\
    \nonumber
    & = &
    -\E [ L (-L)^{\alpha - 1 } ( X - \E [ X ] ) (-L)^{-\alpha} ( Y - \E [ Y ] ) ] 
    \\
  \nonumber
    & = &
     \E [ \nabla^* \nabla (-L)^{\alpha - 1} ( X - \E [ X ] ) (-L)^{-\alpha} ( Y - \E [ Y ] ) ] 
    \\
    \nonumber
    & = &
    \frac{1}{2}
    \E\[\int_0^\infty (\nabla_t(-L)^{\alpha - 1 }X)(\nabla_t(-L)^{-\alpha} Y) dt \]. 
\end{eqnarray} 
\end{Proof}
Although $\nabla_t$ does not satisfy the chain rule of derivation, 
we have the following lemma. 
\begin{lemma}
\label{lem:pr}
 The finite difference operator $\nabla$ satisfies the relation  
\begin{align}\label{eq:prod1}
 \nabla_t ( FG) 
 &= ( F \circ \Phi_t ) \nabla_t G + ( G \circ \Phi_t ) \nabla_t F
 - \frac{1}{2}
\int_{2\lfloor t/2\rfloor}^{2\lfloor t/2\rfloor+2}
 ( \nabla_tF\nabla_tG +  \nabla_uF\nabla_uG  )
du,
\end{align}
$t\in \real_+$, provided that
$( F\circ \Phi_s)_{s\in \real_+}$,
$( G\circ \Phi_s)_{s\in \real_+}$
and
$( F^2 \circ \Phi_s)_{s\in \real_+}$,
$( G^2 \circ \Phi_s)_{s\in \real_+}$
are integrable on $[2 n-2,2 n]$, $n\geq 1$, 
$\P$-a.s. 
\end{lemma}
\begin{Proof}
 By \eqref{eq:nablaaspsi}, we have  
\begin{eqnarray*} 
  \lefteqn{
    \nabla_t ( FG) 
  =  \frac{1}{2}
 \int_{2\lfloor t/2\rfloor}^{2\lfloor t/2\rfloor+2}
 ( (FG) \circ \Phi_t - (FG) \circ \Phi_u ) du
  }
  \\
 & = &
 \frac{1}{2}
 \int_{2\lfloor t/2\rfloor}^{2\lfloor t/2\rfloor+2}
 ( F \circ \Phi_u ) ( G \circ \Phi_t - G \circ \Phi_u )
 du
  + \frac{1}{2}
 \int_{2\lfloor t/2\rfloor}^{2\lfloor t/2\rfloor+2}
 ( G \circ \Phi_t ) ( F \circ \Phi_t - F \circ \Phi_u )
 du
 \\
 & = &
 \frac{1}{2}
 ( F \circ \Phi_t )
 \int_{2\lfloor t/2\rfloor}^{2\lfloor t/2\rfloor+2}
 ( G \circ \Phi_t - G \circ \Phi_u )
 du
  + \frac{1}{2}
 ( G \circ \Phi_t )
 \int_{2\lfloor t/2\rfloor}^{2\lfloor t/2\rfloor+2}
 ( F \circ \Phi_t - F \circ \Phi_u )
 du
 \\
 & &
 - \frac{1}{2}
 \int_{2\lfloor t/2\rfloor}^{2\lfloor t/2\rfloor+2}
 ( F \circ \Phi_t - F \circ \Phi_u )
 ( G \circ \Phi_t - G \circ \Phi_u )
 du
 \\
 & = &
  ( F \circ \Phi_t ) \nabla_t G + ( G \circ \Phi_t ) \nabla_t F
 - \frac{1}{2}
 \int_{2\lfloor t/2\rfloor}^{2\lfloor t/2\rfloor+2}
 ( F \circ \Phi_t - F \circ \Phi_u )
 ( G \circ \Phi_t - G \circ \Phi_u )
 du.
\end{eqnarray*}
 Furthermore, we have
\begin{eqnarray*}
  \lefteqn{
\! \! \! \! \! \! \!   \! \! \! \!   \! \! \! \! \! \! \! \! \! \! \! \! \! \! \! \! \! \! \! \! \! \! \! \! \! \! \! \! \! \! \! \! \! \! \! \! \! \! \! \! \! \! \! \! \! \! \! \! \! \! \! \! \! \!
    \int_{2\lfloor t/2\rfloor}^{2\lfloor t/2\rfloor+2}
 ( F \circ \Phi_t - F \circ \Phi_u )
 ( G \circ \Phi_t - G \circ \Phi_u )
     du
 = \int_{2\lfloor t/2\rfloor}^{2\lfloor t/2\rfloor+2}
 ( \nabla_tF -  \nabla_uF  )
 ( \nabla_tG - \nabla_uG )
 du
  }
  \\
  & = & \int_{2\lfloor t/2\rfloor}^{2\lfloor t/2\rfloor+2}
 ( \nabla_tF\nabla_tG +  \nabla_uF\nabla_uG  )
 du,
\end{eqnarray*}
from the equality $ \int_{2\lfloor t/2\rfloor}^{2\lfloor t/2\rfloor+2}\nabla_uFdu=0$.
\end{Proof}
 The next result is a fourth order moment bound
 stated in terms of the gradient operator $\nabla$. 
\begin{prop}\label{prop:EX^4<}
  For any $X\in L^4(\Omega )$ we have 
 \begin{align}\label{eq:EX^4<}\E\[X^4\]\leq
   36
   \E [ \Vert \nabla X \Vert^4_{L^2(\real_+)} ]
   +15
   \E [ \Vert \nabla X \Vert^4_{L^4(\real_+)} ]
   +2\(\E\[X^2\]\)^2. 
 \end{align}
  \end{prop}
\begin{Proof}
  By the covariance relation
  \eqref{eq:intnablaX}, we have
\begin{align*}
\E\[X^4\]&=\Var \[X^2\]+\(\E\[X^2\]\)^2\\
&=\frac{1}{2} \E\[\int_0^\infty\nabla_t\(X^2-\E\[X^2\]\)\nabla_tL^{-1}\(X^2-\E\[X^2\]\)dt\]+\(\E\[X^2\]\)^2\\
&\leq\frac{1}{2} \sqrt{\E\[\int_0^\infty|\nabla_t(X^2)|^2dt\]\E\[\int_0^\infty \|\nabla_tL^{-1}\(X^2-\E\[X^2\]\)\|^2dt\]}+\(\E\[X^2\]\)^2
\\
&\leq\frac{1}{2} \E\[\int_0^\infty|\nabla_t(X^2)|^2dt\]+\(\E\[X^2\]\)^2, 
\end{align*}
 where we applied \eqref{alh}, \eqref{eq:chaosnorm} and \eqref{L-1}.
 Since 
\begin{align*}
 ( X\circ \Phi_t ) \nabla_t X &=X\nabla_tX+(X\circ \Phi_t-X)\nabla_tX
 \\
 &=X\nabla_tX+
 \bigg(
 \nabla_tX-\bigg(X-\frac12\int_{2\lfloor t/2\rfloor}^{2\lfloor t/2\rfloor+2}X\circ\Phi_udu\bigg)
 \bigg)
 \nabla_tX,
\end{align*}
by the relations \eqref{eq:prod1} and \eqref{ex} 
and the bound $(a+b+c)^2\leq 3\(a^2+b^2+c^2\)$, we have 
\begin{align*}
& 
    E\[\int_0^\infty|\nabla_t(X^2)|^2dt\]
  =\E\[\int_0^\infty\bigg( 2( X\circ \Phi_t ) \nabla_t X 
 -  \frac{1}{2}
\int_{2\lfloor t/2\rfloor}^{2\lfloor t/2\rfloor+2}
\big(
| \nabla_tX|^2+|\nabla_uX |^2
\big)
du\bigg)^2dt\]
  \\
 & \leq 3 \E\Bigg[ 4 \int_0^\infty \(X \nabla_t X\)^2 dt 
+4 \int_0^\infty
\left( 
 \bigg(
 \nabla_tX-\bigg(X-\frac12\int_{2\lfloor t/2\rfloor}^{2\lfloor t/2\rfloor+2}X\circ\Phi_udu\bigg)
 \bigg)
 \nabla_tX
 \right)^2
 dt
\\
& \quad
+
\frac{1}{4}
\int_0^\infty \bigg(
\int_{2\lfloor t/2\rfloor}^{2\lfloor t/2\rfloor+2}
 \big( | \nabla_tX|^2+|\nabla_uX |^2
 \big)
 du\bigg)^2 dt 
 \Bigg]
  \\
 & = 12 \E\Bigg[\int_0^\infty \(X \nabla_t X\)^2 dt 
    \Bigg]
  \\
   & \quad + 12 \E \Bigg[ 
\int_0^\infty
\int_{2\lfloor t/2\rfloor}^{2\lfloor t/2\rfloor+2}
\bigg( 
 \bigg(
 \nabla_tX-\bigg(X\circ \Phi_v -\frac12\int_{2\lfloor t/2\rfloor}^{2\lfloor t/2\rfloor+2}X\circ\Phi_udu\bigg)
 \bigg)
 \nabla_tX
 \bigg)^2
 \frac{dv}{2}
 dt
 \Bigg]
  \\
& \quad
+
\frac{3}{4}
\E \Bigg[
  \int_0^\infty \bigg(
\int_{2\lfloor t/2\rfloor}^{2\lfloor t/2\rfloor+2}
 \big( | \nabla_tX|^2+|\nabla_uX |^2
 \big)
 du\bigg)^2 dt 
 \Bigg]
  \\
 & = 12 \E\Bigg[X^2\int_0^\infty \( \nabla_t X\)^2  dt \Bigg] 
 + 12 \E \Bigg[ \int_{2\lfloor t/2\rfloor}^{{2\lfloor t/2\rfloor}+2}
\big( \(\nabla_tX-\nabla_vX\)\nabla_tX\big)^2\frac{dv}2\,dt
\Bigg]
\\
& \quad
+
\frac{3}{4}
\E \Bigg[
  \int_0^\infty \bigg(
\int_{2\lfloor t/2\rfloor}^{2\lfloor t/2\rfloor+2}
 \big( | \nabla_tX|^2+|\nabla_uX |^2
 \big)
 du\bigg)^2 dt 
\Bigg]
 \\
 &\leq  12\sqrt{\E[X^4\,]\E\[\(\int_0^\infty|\nabla_tX|^2dt\)^2\]}
 +12\E\[\int_0^\infty|\nabla_tX|^2\int_{2\lfloor t/2\rfloor}^{{2\lfloor t/2\rfloor}+2}
 \big( |\nabla_tX|^2 + |\nabla_vX|^2 \big) \frac{dv}2 \,dt\]
 \\
 &  \quad
 +3\E\[\int_0^\infty(\nabla_tX)^4dt\]
 \\
 &\leq  12\sqrt{\E[X^4\,]\E\[\(\int_0^\infty|\nabla_tX|^2dt\)^2\]}+15\E\[\int_0^\infty(\nabla_tX)^4dt\].
\end{align*}
Thus, we get
$$\E\[X^4\]\leq 6\sqrt{\E\[X^4 \]\E\[\(\int_0^\infty|\nabla_tX|^2dt\)^2\]}
  +\frac{15}2\E\[\int_0^\infty( \nabla_tX)^4dt\]+\(\E\[X^2\]\)^2.$$
  Denoting
  $$
  a=\sqrt{\E\[\(\int_0^\infty|\nabla_tX|^2dt\)^2\]},
  \qquad
  b=\frac{15}{2} \E\[\int_0^\infty( \nabla_tX)^4dt\]+\(\E\[X^2\]\)^2
  $$
  and
  $x=\sqrt{\E \[X^4\]}$, we rewrite the last inequality as $x^2\leq 6ax+b$, which gives
  $x\leq 3a + \sqrt{9a^2+b}$ and consequently $x^2\leq 2(9a^2+b)+18a^2=36a^2+2b$,
  which yields \eqref{eq:EX^4<}. 
\end{Proof}

\section{Berry-Esseen bound}
\label{s3}
Our main result is a Berry-Esseen bound on the Kolmogorov distance
$d_K ( X,{\mathcal N})$ between the standard normal
distribution $\mathcal N$ on $\real$
and a general functional $X$ of the
uniform i.i.d. sequence $(U_k)_{k\in \inte}$ 
on $[-1,1]$, using the operators $\nabla$ and $L$. 
This result extends Proposition~4.1 in \cite{reichenbachsAoP}, 
see also Theorem~3.1 in \cite{reichenbachs} and Proposition~2.1 in \cite{PS2}, 
from functionals of Bernoulli sequences to more general functionals of
independent random variables. 
\begin{theorem}
   \label{thm:main}
 Let $X\in \Dom ( \nabla )$ be such that $\E[X]=0$. We have 
\begin{align} 
\label{dwf} 
 & 
  d_K (X,\mathcal N) \leq |1- \E [ X^2 ] | +\sqrt{\Var\left[\int_0^\infty
    \nabla_tX  \nabla_t L^{-1} X  \frac{dt}{2}
    \right]}
\\
\nonumber
&
 +
\frac32\sqrt{\,\E  \int_0^\infty
  (\nabla_tX)^4dt}
 \Bigg( \(\E\[X^4\]\,\E \[ \(\int_0^\infty
 |\nabla_t L^{-1}X |^2dt\)^2 \]\)^{1/4}
 \hskip-0.4cm
 +\frac{\sqrt{\pi}}{2}\sqrt{\E \[ 
( (-L)^{-1/2} X )^2 
  \]}\Bigg)
  \\
  \nonumber
&  
 +4\left( \E\[\int_{0}^\infty\(\( {\rm I}+2(-L)^{1/2}\)\(|\nabla_tX|^2 \)\)^2dt\]\,\E\[\int_{0}^\infty\(\( {\rm I}+2(-L)^{1/2} \)\((\nabla_tL^{-1}X)^2 \)\)^2dt\]\right)^{1/4}
.
\end{align} 
\end{theorem}
\begin{Proof}
  For any $x\in \real$,
  let $f_x$ denote the unique bounded solution of the Stein equation
  \begin{equation}
    \label{steineq}
  f_x'(z)-zf_x(z)=\mathbbm{1}_{\{z\leq x\}}-\mathbb P\(\mathcal N\leq x\),
\end{equation} 
  which is continuous, infinitely differentiable on $\real \setminus \{ x\}$,
  and satisfies $0 < f_x(y) < \sqrt{\pi/8}$ and $|f'_x(y)|\leq 1$, $y \in \real$,
  see Lemmas~2.2 and 2.3 in \cite{chenbk}.
  From the Stein equation \eqref{steineq} we have the bound 
  $$d_K(X,\mathcal N)\leq \sup_{x\in \R} \E [ f_x'(X) - X f_x (X) ].
  $$
 For every $f \in {\cal C}^1(\R )$, 
 the finite difference operator $\nabla$ satisfies 
\begin{eqnarray*} 
 \nabla_t f( X)
 & = & \frac{1}{2}
 \int_{2\lfloor t/2\rfloor}^{2\lfloor t/2\rfloor+2}
 ( f(X \circ \Phi_t) - f(X \circ \Phi_s) ) ds
 \\
 & = & \frac{1}{2}
 \int_{2\lfloor t/2\rfloor}^{2\lfloor t/2\rfloor+2}
 \int_{X \circ \Phi_s - X}^{X \circ \Phi_t - X} f'(X+u) du ds
 \\
 & = & \frac{1}{2}
 \int_{2\lfloor t/2\rfloor}^{2\lfloor t/2\rfloor+2}
 \left(
 \int_{X \circ \Phi_s - X}^{X \circ \Phi_t - X} ( f'(X+u) - f'(X) ) du 
 +
 \int_{X \circ \Phi_s - X}^{X \circ \Phi_t - X} f'(X) du
 \right)
 ds\\
 & = &
 f'(X)\nabla_t X 
 +
 \frac{1}{2}
 \int_{2\lfloor t/2\rfloor}^{2\lfloor t/2\rfloor+2}
 \int_{X \circ \Phi_s - X}^{X \circ \Phi_t - X} ( f'(X+u) - f'(X) ) du\,ds, 
 \quad
 t \in \real_+.
\end{eqnarray*} 
 hence 
  by the duality relation \eqref{eq:duality}, we have
\begin{eqnarray} 
  \nonumber
  \lefteqn{
    \E [ f'(X) - X f (X) ] 
=
\E [ f'(X) - f(X) (-\nabla^*\nabla) L^{-1} X ] 
  }
  \\
\nonumber
  & = & 
  \E \left[ f'(X) - \frac{1}{2} \int_0^\infty
    \nabla_t f (X) ( - \nabla_t L^{-1} X ) dt
    \right] 
  \\
\nonumber
  &= & 
  \E \left[ f'(X)\(1 - \frac{1}{2}\int_0^\infty
    \nabla_tX  ( - \nabla_t L^{-1} X ) dt\)
    \right] 
   \\
\nonumber
& & +\frac{1}{4} 
  \E \bigg[ \int_0^\infty
    \int_{2\lfloor t/2\rfloor}^{2\lfloor t/2\rfloor+2}
 \int_{X \circ \Phi_s - X}^{X \circ \Phi_t - X} ( f'(X+u) - f'(X) ) du 
 ds
 \nabla_t L^{-1} X
 dt
 \bigg] 
. 
\end{eqnarray} 
By the covariance relation \eqref{eq:intnablaX} applied with $\alpha=0$ and
the fact that $\E[X]=0$, we have
$$
\E\[X^2\] = \E\[\int_0^\infty (\nabla_tX) ( - \nabla_tL^{-1}X)\frac{dt}{2}\],
$$
 hence from the bound $\|f_x'\|_\infty\leq1$ and Jensen's inequality we obtain 
\begin{align*}
&\left|\E \left[ f'(X)\bigg(1 - \frac{1}{2}\int_0^\infty
    \nabla_tX  ( - \nabla_t L^{-1} X ) dt\bigg)
    \right]\right|\\
    & \leq \E \left[\left|1 - \frac{1}{2}\int_0^\infty
    \nabla_tX  ( - \nabla_t L^{-1} X ) dt\right|
    \right]\\
  & \leq
  |1- \E [ X^2 ] | +
    \E \left[\left| \frac{1}{2}\int_0^\infty
      \nabla_tX  ( - \nabla_t L^{-1} X ) dt
      -
      \E\[\int_0^\infty (\nabla_tX) ( - \nabla_tL^{-1} X)\frac{dt}{2}\]
      \right|
    \right]
  \\
& \leq |1- \E [ X^2 ] | +\Var\left[\int_0^\infty
    \nabla_tX  ( - \nabla_t L^{-1} X ) \frac{dt}{2}
    \right].
\end{align*}
Next, from the Stein equation \eqref{steineq} we have 
$$ 
  \int_{X \circ \Phi_s - X}^{X \circ \Phi_t - X} ( f'_x(X+u) - f'_x(X) ) du
   = A_{s,t}(x,X) + B_{s,t}(x,X),
   \quad 
   x\in \real,
   $$ 
  where
$$ 
  A_{s,t}(x,X) := \int_{X \circ \Phi_s - X}^{X \circ \Phi_t - X} ( (X+u) f_x(X+u) - Xf_x(X) ) du
  $$
  and
  $$
  B_{s,t}(x,X) : =\int_{X \circ \Phi_s - X}^{X \circ \Phi_t - X}
  \big( {\bf 1}_{\{ X+u \leq x \} } - {\bf 1}_{\{ X \leq x \} } \big) du.
$$ 
Thus, we get
\begin{align}\label{eq:A+B}
\left|\E [ f'(X) - X f (X) ]\right|&\leq   |1- \E [ X^2 ] | +\Var\bigg[\int_0^\infty
    \nabla_tX  ( - \nabla_t L^{-1} X ) \frac{dt}{2}
    \bigg]\\\nonumber
    &\ \ \ +\frac{1}{4}\bigg| 
  \E \bigg[ \int_0^\infty
 \int_{2\lfloor t/2\rfloor}^{2\lfloor t/2\rfloor+2}
 A_{s,t}(x,X)
 ds \nabla_t L^{-1} X dt
    \bigg]\bigg| \\\nonumber
    &\ \ \ +\frac{1}{4} \bigg|
  \E \bigg[ \int_0^\infty
 \int_{2\lfloor t/2\rfloor}^{2\lfloor t/2\rfloor+2}
 B_{s,t}(x,X)
 ds
 \nabla_t L^{-1} X dt
    \bigg] \bigg|.
\end{align}
Using the inequality
 \begin{align*}
   |(u+w)f_x(u+w)-wf_x(w)|\leq \bigg(
   | w | +
   \frac{\sqrt{2\pi}}{4}
   \bigg)|u|,\qquad u,w\in\mathbb R,
\end{align*}
 see Lemma~2.3 in \cite{chenbk}, we estimate
\begin{eqnarray*}
  \left|A_{s,t}(x,X)\right|& \leq &
 \int_{\min ( X \circ \Phi_s - X , X \circ \Phi_t - X)}^{\max ( X \circ \Phi_s - X , X \circ \Phi_t - X)}
 \bigg(
 |X|
 +
  \frac{\sqrt{2\pi}}{4}
\bigg)
 |u| du 
 \\
   & \leq &
  \left(
  \frac{\sqrt{2\pi}}{4}
  + |X|\right)\int_{|X \circ \Phi_s - X|}^{|X \circ \Phi_t - X|}
 |u| du 
 \\
  & = &
  \frac{1}{2} \bigg(
  \frac{\sqrt{2\pi}}{4}
  + |X|\bigg)
 \left(
 \left| X \circ \Phi_s - X \right|^2
 + \left| X \circ \Phi_t - X \right|^2
   \right).
\end{eqnarray*} 
 Then, by the Cauchy-Schwarz inequality we have
\begin{align}\nonumber
 &\bigg|\E \bigg[ \int_0^\infty
 \int_{2\lfloor t/2\rfloor}^{2\lfloor t/2\rfloor+2}
 A_{s,t}(x,X)
 ds
 \nabla_t L^{-1} X dt
    \bigg]\bigg| \\\nonumber
&\leq \frac{1}{2} \E \bigg[ \int_0^\infty
 \int_{2\lfloor t/2\rfloor}^{2\lfloor t/2\rfloor+2}
 \bigg(
 \frac{\sqrt{2\pi}}{4}
 + |X|\bigg)
 \big( \left| X \circ \Phi_t - X \right|^2
   + \left| X \circ \Phi_s - X \right|^2
   \big)
 |  \nabla_t L^{-1} X | ds\,dt
 \bigg] \\
  \nonumber
    &\leq \frac{1}{2} \sqrt{\E \bigg[ \int_0^\infty
 \int_{2\lfloor t/2\rfloor}^{2\lfloor t/2\rfloor+2}
 \left( \left| X \circ \Phi_t - X \right|^2
   + \left| X \circ \Phi_s - X \right|^2
   \right)^2
 ds\,dt
    \bigg] } \\  
  \label{eq:intA}
  &\quad\quad \times\(\sqrt{2\E \left[ \int_0^\infty
  (X\nabla_t L^{-1} X )^2 dt
    \right]}+\frac{\sqrt{\pi}}{2}\sqrt{\E \left[ \int_0^\infty
  (\nabla_t L^{-1} X )^2 dt
    \right]}\).
    \end{align}
Next, by the inequality $(a+b+c)^2\leq 3(a^2+b^2+c^2)$,
formula \eqref{eq:intnablaX} with $\alpha =0$ and the relation \eqref{ex},
 we get 
  \begin{align*}
    &\E \bigg[ \int_0^\infty
 \int_{2\lfloor t/2\rfloor}^{2\lfloor t/2\rfloor+2}
 \left( \left| X \circ \Phi_t - X \right|^2
   + \left| X \circ \Phi_s - X \right|^2
   \right)^2
 ds\,dt
 \bigg]
    \\
    &= \frac{1}{2} \E \bigg[ \int_0^\infty
 \int_{2\lfloor t/2\rfloor}^{2\lfloor t/2\rfloor+2}
 \int_{2\lfloor t/2\rfloor}^{2\lfloor t/2\rfloor+2}
 \left( \left| X \circ \Phi_t - X\circ \Phi_v \right|^2
   + \left| X \circ \Phi_s - X\circ \Phi_v \right|^2
   \right)^2
 dv\,ds\,dt \bigg] 
   \\  
   &= \frac{1}{2} \E \bigg[
     \int_0^\infty
 \int_{2\lfloor t/2\rfloor}^{2\lfloor t/2\rfloor+2}
 \int_{2\lfloor t/2\rfloor}^{2\lfloor t/2\rfloor+2}
\left( \left| \nabla_tX - \nabla_vX \right|^2
   + \left| \nabla_sX - \nabla_vX \right|^2
   \right)^2
   dv\,ds\,dt
   \bigg]
   \\    
   &\leq \E \bigg[
     \int_0^\infty
 \int_{2\lfloor t/2\rfloor}^{2\lfloor t/2\rfloor+2}
 \int_{2\lfloor t/2\rfloor}^{2\lfloor t/2\rfloor+2}
 \left(|\nabla_tX|^2 +|\nabla_sX|^2+2|\nabla_vX|^2
   \right)^2
   dv\,ds\, dt
   \bigg]
   \\
   &\leq3 \E \bigg[
     \int_0^\infty
 \int_{2\lfloor t/2\rfloor}^{2\lfloor t/2\rfloor+2}
 \int_{2\lfloor t/2\rfloor}^{2\lfloor t/2\rfloor+2}
 (\nabla_tX)^4 +(\nabla_sX)^4+4(\nabla_vX)^4
 dv\,ds\, dt
 \bigg]
   \\
   &=72\,\E \left[
     \int_0^\infty
     (\nabla_tX)^4dt
     \right].
\end{align*}
  Furthermore, by
Proposition~\ref{sdajklsad} 
applied to $X$ and $(-L)^{-1} X$ with $\alpha= 1/2$,
   we have
   $$
      \E\[ \int_0^\infty
\big( \nabla_t L^{-1} X \big)^2 dt
\]=2\E \big[ 
\big( (-L)^{-1/2} X \big)^2 
\big]
$$
    and 
    \begin{align*}
    \E \[ \int_0^\infty
  (X\nabla_t L^{-1} X )^2 dt\]\leq \sqrt{\E\[X^4\]\,\E \[ \(\int_0^\infty
  (\nabla_t L^{-1}X  )^2\)^2 dt\]}.
    \end{align*}
    Applying the last three inequalities to \eqref{eq:intA}, we finally obtain
    \begin{align*}
    &\left|\E \bigg[ \int_0^\infty
 \int_{2\lfloor t/2\rfloor}^{2\lfloor t/2\rfloor+2}
 A_{s,t}(x,X)
 ds \nabla_t L^{-1} X dt
    \bigg]\right|\\
   &\leq 6 \sqrt{\,\E  \int_0^\infty
 (\nabla_tX)^4dt}\Bigg(\(\E\[X^4\]\,\E \[ \(\int_0^\infty
  (\nabla_t L^{-1}X  )^2\)^2 dt\]\)^{1/4}+\frac{\sqrt{\pi}}{2} \sqrt{\E \[ 
( (-L)^{-1/2} X )^2 
    \]}\Bigg).
    \end{align*}
    Regarding the last term in \eqref{eq:A+B}, we use   \eqref{ex} and the equivalence  $(  \nabla_t L^{-1} X )\circ\Phi_v =(  \nabla_t L^{-1} X )$, which is
    valid for $2\lfloor t/2\rfloor \leq v < 2\lfloor t/2\rfloor+2$, and get
\begin{align}\nonumber
& \left|
  \E \bigg[ \int_0^\infty
 \int_{2\lfloor t/2\rfloor}^{2\lfloor t/2\rfloor+2}
 B_{s,t}(x,X)
 ds \nabla_t L^{-1} X dt
    \bigg] \right| \\\nonumber
    &= \bigg|
  \E \bigg[ \int_0^\infty
 \int_{2\lfloor t/2\rfloor}^{2\lfloor t/2\rfloor+2}
 \left(
 \int_{X \circ \Phi_s }^{X \circ \Phi_t } ( {\bf 1}_{\{ u \leq x \} } - {\bf 1}_{\{ X \leq x \} } )  du 
 \right)
 ds \nabla_t L^{-1} X dt
    \bigg] \bigg|\\\nonumber
       &=\frac{1}{2} \bigg|
  \E \bigg[ \int_0^\infty\bigg(
 \int_{2\lfloor t/2\rfloor}^{2\lfloor t/2\rfloor+2} \int_{2\lfloor t/2\rfloor}^{2\lfloor t/2\rfloor+2}
 \int_{X \circ \Phi_s }^{X \circ \Phi_t } ( {\bf 1}_{\{ u \leq x \} } - {\bf 1}_{\{ X\circ \Phi_v \leq x \} })  du
 ds\,dv\bigg)\, \nabla_t L^{-1} X dt
 \bigg]\bigg|
  \\
  \label{eq:intB=intK}
    &=\frac{1}{2} \Bigg|
  \E \Bigg[ \sum_{m=0}^\infty
 \int_{2m}^{2m+2}K_m(t,X)
 \nabla_t L^{-1} X dt
    \Bigg] \Bigg|,
\end{align}
where 
$$K_m(t,x,X):=\int_{2m}^{2m+2} \int_{2m}^{2m+2}
\int_{X \circ \Phi_s }^{X \circ \Phi_t }
\big(
    {\bf 1}_{\{ u \leq x \} }-{\bf 1}_{\{ X\circ \Phi_v \leq x \} }
    \big)  du\, ds\,dv,\quad 2m \leq t < 2m+2.
    $$
Next, we rewrite $K_m(t,X)$ as follows 
\begin{align}\nonumber
K_m(t,x,X)&= \int_{2m}^{2m+2}
 \int_{X \circ \Phi_s }^{X \circ \Phi_t } \int_{2m}^{2m+2} \big( {\bf 1}_{\{ X \circ \Phi_t  \leq x \} } -{\bf 1}_{\{ X\circ \Phi_v \leq x \} }\big)  dv\,du\, ds\\\nonumber
 &\ \ \ + \int_{2m}^{2m+2}
 \int_{X \circ \Phi_s }^{X \circ \Phi_t } \int_{2m}^{2m+2} \big(  {\bf 1}_{\{ u \leq x \} }-{\bf 1}_{\{ X \circ \Phi_t  \leq x \} } \big)  dv\,du\, ds\\\nonumber
 &=4\nabla_tX\nabla_t{\bf 1}_{\{ X  \leq x \} }+2\int_{2m}^{2m+2}
 \int_{X \circ \Phi_s }^{X \circ \Phi_t }  \big(  {\bf 1}_{\{ u \leq x \} }-{\bf 1}_{\{ X \circ \Phi_t  \leq x \} } \big) \,du\, ds\\\label{eq:aux4}
  &=-4\nabla_tX\nabla_t{\bf 1}_{\{ X  > x \} }+2\int_{2m}^{2m+2}
 \int_{X \circ \Phi_s }^{X \circ \Phi_t }  \big(  {\bf 1}_{\{ u \leq x \} }-{\bf 1}_{\{ X \circ \Phi_t  \leq x \} } \big) \,du\, ds,
\end{align}
where we used the equality $\nabla_t{\bf 1}_{\{ X  \leq x \} }=-\nabla_t{\bf 1}_{\{ X  > x \} }$.
Next, we consider two cases.
 
\noindent
$(i)$ If $X \circ \Phi_t  > x$, 
we have
\begin{align}\nonumber
  K_m(t,x,X)&=-4\nabla_tX\nabla_t{\bf 1}_{\{ X  > x \} }+2 \int_{2\lfloor t/2\rfloor}^{2\lfloor t/2\rfloor+2}
 \int_{X \circ \Phi_s }^{X \circ \Phi_t }
     {\bf 1}_{\{ u \leq x \} } \,du\, ds
  \\\label{eq:aux20}
&=-4\nabla_tX\nabla_t{\bf 1}_{\{ X  > x \} }+2  \int_{2\lfloor t/2\rfloor}^{2\lfloor t/2\rfloor+2} {\bf 1}_{\{ X \circ \Phi_s \leq x \} } 
  ( x - X \circ \Phi_s ) \, ds.
\end{align}
Note that the last expression depends  only
on $m:=\lfloor t/2\rfloor$ and may be bounded for $x< \int_{2m}^{2m+2}X \circ \Phi_u\,du / 2$ as follows
\begin{align*}
  0&\leq \int_{2\lfloor t/2\rfloor}^{2\lfloor t/2\rfloor+2} {\bf 1}_{\{ X \circ \Phi_s \leq x \} } 
  ( x - X \circ \Phi_s ) \, ds
  \\
&= x\int_{2m}^{2m+2} {\bf 1}_{\{ X \circ \Phi_s \leq x \} } \, ds
- \int_{2m}^{2m+2}
   X\circ \Phi_udu+
 \int_{2m}^{2m+2}
 {\bf 1}_{\{ X \circ \Phi_u  > x \} }
   X\circ \Phi_udu\\
   &= \(x- \frac{1}{2} \int_{2m}^{2m+2}X \circ \Phi_u\,du\)\int_{2m}^{2m+2} {\bf 1}_{\{ X \circ \Phi_s \leq x \} } \, ds\\
   &\quad 
 + \int_{2m}^{2m+2}
 {\bf 1}_{\{ X \circ \Phi_u  > x \} }
 X\circ \Phi_udu
    - \frac{1}{2} 
 \int_{2m}^{2m+2}
 {\bf 1}_{\{ X \circ \Phi_s  > x \} }
 ds
 \int_{2m}^{2m+2} X\circ \Phi_udu
 \\
    &\leq 
 \int_{2m}^{2m+2}
 \({\bf 1}_{\{ X \circ \Phi_u  > x \} }
 -\frac{1}{2} \int_{2m}^{2m+2}
 {\bf 1}_{\{ X \circ \Phi_s  > x \} }
   ds\)
   X\circ \Phi_u du
   \\
   &=
 \int_{2m}^{2m+2}\nabla_u{\bf 1}_{\{ X   > x \}}X\circ \Phi_udu\\
  &=
 \int_{2m}^{2m+2}\nabla_u{\bf 1}_{\{ X  > x \}}\nabla_uXdu.
\end{align*}
Consequently, for $x< \int_{2m}^{2m+2}X \circ \Phi_u\,du / 2$ we get
\begin{align}\nonumber
&
 \int_{2m}^{2m+2}{\bf 1}_{\{ X\circ\Phi_t  > x \}}K_m(t,X)
 \nabla_t L^{-1} X dt
   \\\label{eq:intK<}
       &\leq 4\,
 \int_{2m}^{2m+2}\left|\nabla_tX\nabla_t{\bf 1}_{\{ X  > x \} }
 \nabla_t L^{-1} X \right| dt
    \\\nonumber
    &\ \ \ + 2\,
\left| \int_{2m}^{2m+2}\nabla_uX\nabla_u{\bf 1}_{\{ X  > x \}}du\right|\left|\int_{2m}^{2m+2}
 \nabla_t{\bf 1}_{\{ X> x \}} \nabla_t L^{-1} X dt\right|
    ,
\end{align}
where we also changed ${\bf 1}_{\{ X\circ\Phi_t  > x \}}$ into  $\nabla_t{\bf 1}_{\{ X> x \}}$ in the last integral, which is justified by $\int_{2m}^{2m+2}
  \nabla_t L^{-1} X dt=0$.
In order to obtain the same bound in the case $x\geq  \int_{2m}^{2m+2}X \circ \Phi_u\,du / 2$, we rewrite \eqref{eq:aux20} as 
\begin{align*}
K_m(t,x,X)&=-4\nabla_tX\nabla_t{\bf 1}_{\{ X  > x \} }+2  \int_{2\lfloor t/2\rfloor}^{2\lfloor t/2\rfloor+2} {\bf 1}_{\{ X \circ \Phi_s \leq x \} } 
  ( x - X \circ \Phi_t+\nabla_tX-\nabla_sX ) \, ds\\
  &=2  \int_{2\lfloor t/2\rfloor}^{2\lfloor t/2\rfloor+2} {\bf 1}_{\{ X \circ \Phi_s \leq x \} } 
  ( x - X \circ \Phi_t-\nabla_sX ) \, ds\\
&=2
 \int_{2m}^{2m+2}
 \nabla_s{\bf 1}_{\{ X   > x \} }
   \nabla_sX 
 ds-2\int_{2m}^{2m+2} {\bf 1}_{\{ X \circ \Phi_s \leq x \} } 
 \({X \circ \Phi_t }-x\)  \, ds,
\end{align*}
 and we estimate the last integral by 
\begin{align*}
0&\leq \int_{2m}^{2m+2} {\bf 1}_{\{ X \circ \Phi_s \leq x \} } 
\({X \circ \Phi_t }-x\) \, ds\\
 &\leq  \int_{2m}^{2m+2} {\bf 1}_{\{ X \circ \Phi_s \leq x \} } 
 \, ds\(X \circ \Phi_t-\frac{1}{2} \int_{2m}^{2m+2}X \circ \Phi_u\,du\)\\
 &=-\nabla_tX\nabla_t{\bf 1}_{\{ X  \leq x \} } =\nabla_tX\nabla_t{\bf 1}_{\{ X > x \} },
\end{align*}
which shows that the inequality \eqref{eq:intK<} is valid for all $x\in \R$
under the condition $x<X \circ \Phi_t$.
Thus, applying the Cauchy-Schwarz inequality several times 
and using  the bound $|\nabla_t{\bf 1}_{\{ X  \leq x \} }|\leq1$,
we obtain
\begin{align*}\nonumber
&\left|\E \left[ \sum_{m=0}^\infty
 \int_{2m}^{2m+2}{\bf 1}_{\{ X\circ\Phi_t  > x \}}K_m(t,x,X)
 \nabla_t L^{-1} X dt
    \right]\right|\\
       &\leq4\sqrt{\E \left[ \sum_{m=0}^\infty
 \int_{2m}^{2m+2}\left|\nabla_u{\bf 1}_{\{ X   > x \}}\right||\nabla_uX|^2 du
    \right]\,\E \left[ \sum_{m=0}^\infty
 \int_{2m}^{2m+2}\left|\nabla_u{\bf 1}_{\{ X   > x \}}\right|(\nabla_uL^{-1}X)^2du
    \right]}\\\nonumber
    &\ \ \ +4 \sqrt{\E \left[ \sum_{m=0}^\infty
 \int_{2m}^{2m+2}(\nabla_u{\bf 1}_{\{ X   > x \}})^2|\nabla_uX|^2du
    \right]\,\E \left[ \sum_{m=0}^\infty
 \int_{2m}^{2m+2}(\nabla_u{\bf 1}_{\{ X   > x \}})^2(\nabla_uL^{-1}X)^2du
    \right]}\\
    &\leq8\sqrt{\E \left[ 
 \int_{0}^{\infty}\left|\nabla_u{\bf 1}_{\{ X   > x \}}\right||\nabla_uX|^2 du
    \right]\,\E \left[
 \int_{0}^{\infty}\left|\nabla_u{\bf 1}_{\{ X   > x \}}\right|(\nabla_uL^{-1}X)^2du
    \right]}.
\end{align*}
By the duality relation \eqref{eq:duality},
H\"older's inequality and the formula \eqref{eq:Enabla*^2}, we get
\begin{align}\nonumber
&\E \left[ 
 \int_{0}^{\infty}\left|\nabla_u{\bf 1}_{\{ X   > x \}}\right||\nabla_uX|^2 du
    \right]\\\nonumber
    &=\E \left[ 
 \int_{0}^{\infty}\nabla_u{\bf 1}_{\{ X   > x \}}\text{sgn}({\nabla_u{\bf 1}_{\{ X   > x \}}})|\nabla_uX|^2 du
    \right]\\\nonumber
  &=2\E \left[ 
 {\bf 1}_{\{ X   > x \}}\nabla^*\(\text{sgn}({\nabla_u{\bf 1}_{\{ X   > x \}}})|\nabla_uX|^2 \)
    \right]\\\nonumber
      &\leq2\sqrt{\E \Big[ 
      \big(
      \nabla^*\big(\text{sgn}({\nabla_u{\bf 1}_{\{ X   > x \}}})|\nabla_uX|^2 \big)
      \big)^2
    \Big]}\\\label{eq:aux9}
    &=2\sqrt{\E\[\int_0^\infty (\nabla_uX)^4dt\]+\E\[\int_{0}^\infty\int_{0}^\infty\(\nabla_s\(\text{sgn}({\nabla_u{\bf 1}_{\{ X   > x \}}})|\nabla_uX|^2 \)\)^2ds\,du\]}. 
\end{align}
Next, we observe that by the covariance relation \eqref{eq:intnablaX} with $\alpha=\frac12$,
we have 
\begin{align*}
&\E\[\int_{0}^\infty\int_{0}^\infty\(\nabla_s\(\text{sgn}({\nabla_u{\bf 1}_{\{ X   > x \}}})|\nabla_uX|^2 \)\)^2ds\,du\]\\
&=\E\[\int_{0}^\infty {\bf 1}_{\{{\nabla_u{\bf 1}_{\{ X   > x \}}}>0\}}\int_{0}^\infty\(\nabla_s\(|\nabla_uX|^2 \)\)^2ds+{\bf 1}_{\{{\nabla_u{\bf 1}_{\{ X   > x \}}}<0\}}\int_{0}^\infty\(\nabla_s\(-|\nabla_uX|^2 \)\)^2ds\,du\]\\
  &\leq \E\[\int_{0}^\infty\int_{0}^\infty\(\nabla_s\(|\nabla_uX|^2 \)\)^2ds\,du\]
  \\
   & =2 \E\[\int_{0}^\infty\((-L)^{1/2}\(|\nabla_uX|^2 \)\)^2du\].
\end{align*}
Applying this to \eqref{eq:aux9}, we get
$$\E \left[ 
 \int_{0}^{\infty}\left|\nabla_u{\bf 1}_{\{ X   > x \}}\right||\nabla_uX|^2 du
 \right]\leq 2\sqrt{ \,\E\[\int_{0}^\infty\(\({\rm I}+2(-L)^{1/2}\)\(|\nabla_uX|^2 \)\)^2du\]},
$$
 and analogously we obtain
$$\E \left[ 
 \int_{0}^{\infty}\left|\nabla_u{\bf 1}_{\{ X   > x \}}\right|(\nabla_uL^{-1}X)^2 du
    \right]\leq 2\sqrt{ \,\E\[\int_{0}^\infty\(\({\rm I}+2(-L)^{1/2}\)\((\nabla_uL^{-1}X)^2 \)\)^2du\]},$$
    which eventually gives us
    \begin{align}\label{eq:aux6}
     &\left|\E \left[ \sum_{m=0}^\infty
 \int_{2m}^{2m+2}{\bf 1}_{\{ X\circ\Phi_t  > x \}}K_m(t,x,X)
 \nabla_t L^{-1} X dt
    \right]\right|\\\nonumber
    &\leq  16\Bigg( \E\[\int_{0}^\infty\(\({\rm I}+2(-L)^{1/2}\)\(|\nabla_uX|^2 \)\)^2du\]\\
    &\ \ \ \ \ \ \ \ \ \ \ \ \times\E\[\int_{0}^\infty\(\({\rm I}+2(-L)^{1/2}\)^{1/2}\((\nabla_uL^{-1}X)^2 \)\)^2du\]\Bigg)^{1/4}.
    \end{align}
    \noindent
$(ii)$ In case $X \circ \Phi_t  \leq x$ we observe that, denoting
  $$\widetilde{K}_m(t,x,X) :=-4\nabla_tX\nabla_t{\bf 1}_{\{ X  \geq x \} }+2\int_{2m}^{2m+2}
    \int_{X \circ \Phi_s }^{X \circ \Phi_t }  \big(
        {\bf 1}_{\{ u < x \} }-{\bf 1}_{\{ X \circ \Phi_t  < x \} } \big) \,du\, ds,$$
 which comes from  \eqref{eq:aux4} by changing weak
 inequalities into strict ones and conversely, and repeating all the above argument, we arrive at
    \begin{align}\label{eq:aux7}
     &\left|\E \left[ \sum_{m=0}^\infty
 \int_{2m}^{2m+2}{\bf 1}_{\{ X\circ\Phi_t  \geq x \}}\widetilde{K}_m(t,x,X)
 \nabla_t L^{-1} X dt
    \right]\right|\\\nonumber
    &\leq  16\( \E\[\int_{0}^\infty\(\({\rm I}+2(-L)^{1/2}\)\(|\nabla_uX|^2 \)\)^2du\]\,\E\[\int_{0}^\infty\(\({\rm I}+2(-L)^{1/2}\)\((\nabla_uL^{-1}X)^2 \)\)^2du\]\)^{1/4}.
    \end{align}
    Next, by \eqref{eq:aux4} we have,
     for $m=\lfloor t/2\rfloor$ and $X\circ\Phi_t  \leq x$, 
\begin{align*}
K_m(t,x,X) & =-4\nabla_tX\nabla_t{\bf 1}_{\{ X  > x \} }+2\int_{2m}^{2m+2}
  \int_{X \circ \Phi_s }^{X \circ \Phi_t }
  \big( {\bf 1}_{\{ u \leq x \} }-{\bf 1}_{\{ X \circ \Phi_t  \leq x \} } \big) \,du\, ds\\
 &=4\nabla_tX\nabla_t{\bf 1}_{\{ X  \leq x \} }-\int_{2m}^{2m+2}
 \int_{X \circ \Phi_s }^{X \circ \Phi_t }  \big(  {\bf 1}_{\{ u \geq x \} }-{\bf 1}_{\{ X \circ \Phi_t  \geq x \} } \big) \,du\, ds\\
  &=4\nabla_tX\nabla_t{\bf 1}_{\{ -X  \geq -x \} }-\int_{2m}^{2m+2}
 \int_{X \circ \Phi_s }^{X \circ \Phi_t }  \big(  {\bf 1}_{\{ -u \leq - x \} }-{\bf 1}_{\{ -X \circ \Phi_t  \leq - x \} } \big) \,du\, ds\\
   &=-4\nabla_t(-X)\nabla_t{\bf 1}_{\{ -X  \geq -x \} }+\int_{2m}^{2m+2}
 \int_{-X \circ \Phi_s }^{-X \circ \Phi_t }  \big( {\bf 1}_{\{ u \leq - x \} }-{\bf 1}_{\{ -X \circ \Phi_t  \leq - x \} } \big) \,du\, ds\\
 &=\widetilde{K}_m(t,-x,-X).
\end{align*}
Thus, using \eqref{eq:aux7} with $-x$ and $-X$ instead of $x$ and $X$
respectively, we get
 \begin{align}\nonumber
     &\left|\E \left[ \sum_{m=0}^\infty
 \int_{2m}^{2m+2}{\bf 1}_{\{ X\circ\Phi_t  \leq x \}} K_m(t,x,X)
 \nabla_t L^{-1} X dt
    \right]\right|\\\nonumber
    &=\left|\E \left[ \sum_{m=0}^\infty
 \int_{2m}^{2m+2}{\bf 1}_{\{ -X\circ\Phi_t  \geq -x \}} \widetilde{K}_m(t,-x,-X)
 \nabla_t L^{-1} (-X) dt
    \right]\right|\\\label{eq:aux8}
    &\leq 16\( \E\[\int_{0}^\infty\(\({\rm I}+2(-L)^{1/2}\)\(|\nabla_uX|^2 \)\)^2du\]\,\E\[\int_{0}^\infty\(\({\rm I}+2(-L)^{1/2}\)\((\nabla_uL^{-1}X)^2 \)\)^2du\]\)^{1/4}.
    \end{align}
Combining \eqref{eq:aux6} and \eqref{eq:aux8} with \eqref{eq:intB=intK}, we finally obtain
    \begin{align*}
&\frac{1}{4} \left|
  \E \left[ \int_0^\infty
 \int_{2\lfloor t/2\rfloor}^{2\lfloor t/2\rfloor+2}
 B_{s,t}(x,X)
 ds
 \nabla_t L^{-1} X dt
    \right] \right| \\
    &=\frac{1}{8}  \left|
  \E \left[ \sum_{m=0}^\infty
 \int_{2m}^{2m+2}\({\bf 1}_{\{ X\circ\Phi_t  > x \}}+{\bf 1}_{\{ X\circ\Phi_t  \leq x \}}\)K_m(t,X)
 \nabla_t L^{-1} X dt
    \right] \right|\\
    &\leq 4\( \E\[\int_{0}^\infty\(\({\rm I}+2(-L)^{1/2}\)\(|\nabla_uX|^2 \)\)^2du\]\,\E\[\int_{0}^\infty\(\({\rm I}+2(-L)^{1/2}\)\((\nabla_uL^{-1}X)^2 \)\)^2du\]\)^{1/4}, 
\end{align*}
   which ends the proof.
    \end{Proof}
\section{Sums of multiple stochastic integrals} 
\label{s4}
The next proposition applies Theorem~\ref{thm:main}
 to sums of multiple stochastic integrals.
 It extends Theorem~3.1 of \cite{PS2}
from functionals of Bernoulli sequences to functionals of
independent random variables, 
see also earlier results such as Proposition~3.7 in \cite{nourdinpeccati}
in the case of multiple Wiener integrals. 
\begin{prop}
  \label{prop:dK<Sum}
  For any $X\in L^2 (\Omega )$ written as a sum $X=\sum_{k=1}^d I_k(f_k)$
  of multiple stochastic integrals where  
  $f_k\in \widehat{L}^2(\R_+^k)$ satisfies \eqref{ass:int0},
  $k=1,\ldots ,d$,
  we have 
\begin{eqnarray*}
  \lefteqn{
    d_K(X , \mathcal{N})  \leq   
 \E\[\left|1 - \E [ X^2 ]     \right|     \]
  }
  \\
 & &
 +
  C_d \sqrt{
    \sum_{0\leq l< i\leq d}\left\| f_{i}\star_{i}^l f_{i}\right\|^2_{
      L^2(\R_+^{i-l})}+\sum_{1\leq l< i\leq d}
    \(\left\| f_{i}\star_{l}^l f_{i}\right\|^2_{
      L^2(\R_+^{2(i-l)})}+\left\| f_{l}\star_{l}^l f_{i}\right\|^2_{
      L^2(\R_+^{i-l})}\)}, 
\end{eqnarray*}
for some $C_d>0$. 
\end{prop}
\begin{Proof}
    Since $|\nabla_t X|^2$ and $(\nabla_t L^{-1}X)^2$ are sums of multiple integrals
  of orders $2d-2$ and below,
  the relation~\eqref{eq:chaosnorm} shows the bound 
  $$\E\big[ \(\({\rm I}+2(-L)^{1/2}\)\(|\nabla_tX|^2 \)\)^2\big]
  \leq 2d\E\[\(\nabla_tX\)^4\],$$
    and
    $$\E\big[\(\({\rm I}+2(-L)^{1/2}\)\((\nabla_tL^{-1}X)^2 \)\)^2\big]\leq 2d\E\big[\(\nabla_tL^{-1}X\)^4\big].$$ 
    Additionally, by \eqref{eq:chaosnorm} we also have
    $$\E \[ ( (-L)^{-1/2} X )^2\]\leq \E \[X ^2\]\leq \sqrt{\E\[X^4\]}.$$
    Applying these  inequalities to \eqref{dwf} in Theorem~\ref{thm:main}, we get
\begin{align*} 
d_K (X,\mathcal N)
 &\leq  |1- \E [ X^2 ] | +\sqrt{\Var\left[\int_0^\infty
    \nabla_tX  \nabla_t L^{-1} X {dt}
    \right]}
\\
&
 \ \ \ +
 \frac32
 \(\E\[X^4\]\)^{1/4}
 \sqrt{\E  \int_0^\infty
  (\nabla_tX)^4dt}
\Bigg(
1 +
\(
\E \[ \(\int_0^\infty
  (\nabla_t L^{-1}X  )^2dt\)^2 \]\)^{1/4} \Bigg)\\
  &\ \ \ +6d\sqrt{\E\[\int_0^\infty\(\nabla_tX\)^4dt\]\E\[\int_0^\infty\(\nabla_tL^{-1}X\)^4dt\]}.
\end{align*}   
 Denoting 
 $$R_X :=\sum_{1\leq i\leq j\leq d}\sum_{k=1}^{i}\sum_{l=0}^k\mathbbm{1}_{\{i=j=k=l\}^c}\left\| f_{i}\star_{k}^l f_{j}\right\|^2_{\widehat{L}^2(\R_+^{i+j-k-l})},
 $$ 
 it follows from the proof of Corollary~3.2 in \cite{PS3} that  
     \begin{align}\label{eq:R_X<}
 R_X\leq c_d \(
    \sum_{0\leq l< i\leq d}\left\| f_{i}\star_{i}^l f_{i}\right\|^2_{
      L^2(\R_+^{i-l})}+\sum_{1\leq l< i\leq d}
    \(\left\| f_{i}\star_{l}^l f_{i}\right\|^2_{
      L^2(\R_+^{2(i-l)})}+\left\| f_{l}\star_{l}^l f_{i}\right\|^2_{
      L^2(\R_+^{i-l})}\)\),
      \end{align}
      and
    \begin{align}\label{eq:...<R_X}
    \Var\left[\int_0^\infty
    \nabla_tX  \nabla_t L^{-1} X {dt}
    \right]\leq c_dR_X,\qquad
       \E  \[\int_0^\infty
 (\nabla_tX)^4dt\]\leq c_d R_X,
 \end{align}
    for some $c_d\geq0$. Taking $L^{-1}X$ as $X$ in the last inequality, we also have
    $$\E  \[\int_0^\infty
 (\nabla_tL^{-1}X)^4dt\]\leq c'_d R_X,$$
for some $c'_d\geq0$.    Furthermore, since
    $$\nabla_tL^{-1}X = \sum_{k=0}^{d-1} I_k \(f_{k+1}(t,\cdot)\)$$
and the functions $f_k$ satisfy \eqref{ass:int0},
 the multiplication formula \eqref{eq:mf} gives  
     \begin{align*}
  \int_0^\infty (\nabla_tL^{-1} X)^2dt&=\int_0^\infty\sum_{0\leq i\leq j < d-1 }\sum_{k=0}^{i}\sum_{l=0}^kc_{i,j,l,k}I_{i+j-k-l}
   \big(
   f_{i+1}(t,\cdot)
   \hskip0.1cm \widetilde{\star}_{k}^l f_{j+1}(t,\cdot)\big)dt\\
   &=\sum_{0\leq i\leq j < d-1 }\sum_{k=0}^{i}\sum_{l=0}^kc_{i,j,l,k}I_{i+j-k-l}
   \(\int_0^\infty
   f_{i+1}(t,\cdot)
   \hskip0.1cm \widetilde{\star}_{k}^l f_{j+1}(t,\cdot)dt\)
   \end{align*} 
    for some $c_{i,j,l,k}\geq 0$, and consequently
      \begin{align*}
        &\E\[\(\int_0^\infty (\nabla_tL^{-1} X)^2dt\)^2\]
        \leq c_d\sum_{0\leq i\leq j < d }\sum_{k=0}^{i}\sum_{l=0}^k\left\|
   \(\int_0^\infty
   f_{i+1}(t,\cdot)
   \hskip0.1cm {\star}_{k}^l f_{j+1}(t,\cdot)dt\)\right\|_{L^2(\R_+^{i+j-k-l})}^2\\
   &= c_d\sum_{0\leq i\leq j < d }\sum_{k=0}^{i}\sum_{l=0}^k\left\|
   \(
   f_{i+1}
   \hskip0.1cm {\star}_{k+1}^{l+1} f_{j+1}\)\right\|_{L^2(\R_+^{i+j-k-l})}^2\\
   &= c_d\(\sum_{1\leq i\leq j < d }\sum_{k=1}^{i}\sum_{1=0}^k\mathbbm{1}_{\{i=j=k=l\}^c}\left\|
   \(
   f_{i}
   \hskip0.1cm {\star}_{k}^{l} f_{j}\)\right\|_{L^2(\R_+^{i+j-k-l})}^2+\sum_{i=1}^d
   \(
   f_{i}
   \hskip0.1cm {\star}_{i}^{i} f_{i}\)^2\)\\
   &\leq c_d \(R_X+(\E [ X^2 ] )^2\).
   \end{align*} 
   Similarly, we get for some $C_{i,j,k,l}\geq0$ 
   \begin{align*}
   \E\[X^4\]&\leq c_d\E\[\(\sum_{0\leq i\leq j < d }\sum_{k=0}^{i}\sum_{l=0}^kC_{i,j,l,k}I_{i+j-k-l}
   \( f_{i}\hskip0.1cm \widetilde{\star}_{k}^l f_{j}\)\)^2\]\\
   &\leq c_d\sum_{0\leq i\leq j < d }\sum_{k=0}^{i}\sum_{l=0}^k
   \left\| f_{i}\hskip0.1cm {\star}_{k}^l f_{j}\right\|^2_{L^2(\R_+^{i+j-k-l})}\\
   &=c_d\(R_X+\sum_{i=1}^d
   \(f_{i}  \hskip0.1cm {\star}_{i}^{i} f_{i}\)^2+\sum_{1\leq i\leq j\leq d}
   \left\| f_{i}\hskip0.1cm {\star}_{0}^0 f_{j}\right\|^2_{L^2(\R_+^{i+j})}\)\\
    &=c_d\(R_X+\sum_{i=1}^d
   \left\|f_i\right\|^4_{L^2(\R^i)}+\sum_{1\leq i\leq j\leq d}
  \left\|f_i\right\|^2_{L^2(\R^i)}\left\|f_j\right\|^2_{L^2(\R^j)}\)\\
  &\leq c_d\(R_X+( \E [ X^2 ] )^2\).
   \end{align*}
   This finally gives us
   $$d_K(X,\mathcal N)\leq |1- \E [ X^2 ] |+c_d\sqrt{R_X}\Big(1+\((R_X+ \ E [ X^2 ] )^{1/4}+1\)^2\Big).
   $$
   Since $d_K(X,\mathcal N)\leq1$, we may assume that $\E [ X^2 ] $
   and $R_X$ are bounded, which implies
    $$d_K(X,\mathcal N)\leq |1- \E [ X^2 ] |+c_d\sqrt{R_X},$$
    and the assertion of the corollary follows from \eqref{eq:R_X<}.
\end{Proof}
\noindent 
Next, due to the identity $\nabla_tL^{-1}I_d(f_d)=I_{d-1}\(f_d(t,*)\)$, $d\geq1$, the bound in Theorem \ref{thm:main}
can be significantly simplified in
the case of multiple stochastic integrals $I_d(f_d)$.
\begin{prop}\label{prop:dKI}
For $X=I_d(f_d)$ a multiple stochastic integral of order $d\geq 1$, we have
 \begin{eqnarray*} 
   \lefteqn{
     \! \! \! \! \! \! 
     d_K (X,\mathcal N)}
   \\
 &\leq &  |1-\E [ X^2 ] | +\frac{1}{d} \sqrt{\Var\left[\int_0^\infty
    (\nabla_tX)^2  \frac{dt}{2}
    \right]}+
\frac{12 +
5\sqrt[4]{\E\[X^4\]}}{\sqrt d}\sqrt{\,\E \left[ \int_0^\infty
    (\nabla_tX)^4dt \right]
}\\
&\leq & |1-\E [ X^2 ] | +\sqrt{\Var\left[\int_0^\infty
    (\nabla_tX)^2  \frac{dt}{2}
    \right]}+24\sqrt{\E \left[ \int_0^\infty
 (\nabla_tX)^4dt \right]}.
\end{eqnarray*} 
\end{prop}
\begin{Proof}
  In view of the relation
  $$
  (-L)^{1/2}I_d(f_d)=\frac{1}{\sqrt d}I_d(f_d)
  $$
  and the covariance identity \eqref{eq:intnablaX} applied with $\alpha=0$, Theorem \ref{thm:main} gives 
\begin{align*} 
& d_K (X,\mathcal N)
\leq  |1-\E [ X^2 ] | +\frac{1}{d} \sqrt{\Var\left[\int_0^\infty
    |\nabla_tX|^2  \frac{dt}{2}
    \right]}
\\
&
+
\frac{3}{2\sqrt d}
\sqrt{\,\E  \int_0^\infty
  (\nabla_tX)^4dt}
\(\(\E\[X^4\]\(\frac{1}{d^2}\Var \[ \int_0^\infty
  |\nabla_t X |^2dt \]+4(\E [ X^2 ] )^2\)\)^{1/4}+\frac{\sqrt{\pi }}{2}\sqrt{\E \[ 
 X ^2 
    \]}\)\\
& +\frac{4}{d} \( \E\[\int_{0}^\infty\(\({\rm I}+2(-L)^{1/2}\)\(|\nabla_tX|^2 \)\)^2dt\]\)^{1/2}
.
\end{align*} 
Since $d_K(X,\mathcal N)\leq1$ by definition,
we may assume that $\sqrt{\Var\left[\int_0^\infty
    |\nabla_tX|^2 dt / 2 
    \right]}\leq d$ and $\E \[ 
 X ^2 
    \]\leq2$. Hence we get
    \begin{align*}
      &
      \! \! \! \! \! \! \! \! \! \! \! \! \!
      \! \! \! \! \! \! \! \! \! \! \! \! \!
      \(\E\[X^4\]\(\frac{1}{d^2} \Var \[ \int_0^\infty
  |\nabla_t X|^2dt \]+4(\E [ X^2 ] )^2\)\)^{1/4}+\frac{\sqrt{\pi }}{2}\sqrt{\E \[ 
 X ^2 
    \]}\\
    &\leq \sqrt[4]{\E\[X^4\]}\(\sqrt[4]{18}+\frac{\sqrt{\pi }}{2}\)\leq \frac{10}{3}\sqrt[4]{\E\[X^4\]}.
    \end{align*}
    Furthermore, since $|\nabla_uX|^2$ is a sum of multiple integrals of orders $2d-2$
    and below, we have by \eqref{eq:chaosnorm} 
    $$
    \E\big[ \big( \big( 2(-L)^{1/2}+{\rm I}\big) (|\nabla_tX|^2 )\big)^2\big]\leq
    \(2\sqrt{2d-2}+1 \)^2
    \E\[\(\nabla_tX\)^4\]\leq 9d\,\E\[\(\nabla_tX\)^4\].$$
    Combining all together we obtain the first inequality from the assertion.
    Next, applying Proposition \ref{prop:EX^4<}, we get
  \begin{align*} 
d_K (X,\mathcal N)
  \leq&  |1-\E [ X^2 ] | +\sqrt{\Var\left[\int_0^\infty
    |\nabla_tX|^2  \frac{dt}{2}
    \right]}
  +  
  \sqrt{\E \left[ \int_0^\infty
  (\nabla_tX)^4dt \right]}\\
  &
\times\Bigg(
12 +5
\Bigg(\({36}\E\[\(\int_0^\infty|\nabla_tX|^2dt\)^2\]+15\E\[\int_0^\infty( \nabla_tX)^4dt\]\)+8\Bigg)^{1/4}\Bigg).
\end{align*} 
Using once again the inequality $d_K(X,\mathcal N)\leq1$, we may assume $\sqrt{\Var\left[\int_0^\infty
    |\nabla_tX|^2 dt / 2 
    \right]}\leq1$ and $ \sqrt{\E \left[ \int_0^\infty
  (\nabla_tX)^4dt \right]}\leq\frac1{17}$, which lets us bound the expressiof in the last parenthesis as follows
  \begin{align*}
  &
12 +
5\Bigg(\({36}\E\[\(\int_0^\infty|\nabla_tX|^2dt\)^2\]+15\E\[\int_0^\infty( \nabla_tX)^4dt\]\)+8\Bigg)^{1/4}\leq 12+5\(81\)^{1/4}=27,
  \end{align*}
  which ends the proof.
\end{Proof}
\section{Applications to $U$-statistics}   
\label{s5}

\subsection{Hoeffding decompositions}
Recall that given $(X_1,\ldots , X_n)$
 a family of independent random variables
and $[n] : = \{1,\ldots ,n\}$, $n\geq 1$, the family
$(\mathcal F_J)_{J\subset [n]}$ of $\sigma$-algebras is defined as 
$$\mathcal F_J:=\sigma (X_j\ : \ j\in J), \qquad J\subset [n]. 
$$
 
\begin{definition} 
  \label{d1}
  A centered $\mathcal F_{[n]}$-measurable
  random variable $W_n$ 
   admits a Hoeffding decomposition if
 it can be written as 
 \begin{align}
   \label{eq:Hd}
 W_n=\sum_{J\subset [n]}W_J,
 \end{align}
 where $(W_J)_{J\subset [n]}$
 is a family of random variables such that 
 $W_J$ is ${\cal F}_J$-measurable, $J\subset [n]$, and 
 $$
 \E\[W_J \mid \mathcal F_K\]=0,
 \qquad
 J\not\subseteq K\subset [n].
 $$
\end{definition}
 For $J=\{k_1,\ldots ,k_{|J|}\} $ with $k_1<k_2<\cdots <k_{|J|}$,
 any $W_J$ in Definition~\ref{d1} can be written as 
 a function $W_J=g_J\big(X_{k_1},\ldots ,X_{k_{|J|}}\big)$
 of $\big(X_{k_1},\ldots ,X_{k_{|J|}}\big)$, with in particular 
  \begin{align}\label{ass:int0W}
    \E\[g_J\(X_j:j\in J\) \mid J\backslash \{k\}\]=0,
    \qquad k\in J, 
 \end{align} 
 and 
 \begin{align}
   \label{eq:Hd2}
   W_n=\sum_{J\subset [n]}
   g_J\big(X_{k_1},\ldots ,X_{k_{|J|}}\big). 
 \end{align}
  Note that if $X_i=U_i$, $i\in[n]$,
  then the chaos decomposition \eqref{eq:decomp}
  coincides with the Hoeffding decomposition \eqref{eq:Hd}, 
  by taking 
  $$W_J:=
  \frac{1}{|J|!}f_{|J|}
  \big(
  2k_1 +1+U_1 ,\ldots ,2k_{|J|-1}+1+U_{|J|-1}, 2k_{|J|}+1+U_{|J|} \big),
  \quad
  J \subset [n], 
  $$
   and Condition~\eqref{ass:int0W} is equivalent to \eqref{ass:int0}.
   
   \medskip
   
  The next Theorem~\ref{thm:S} is a consequence of Proposition~\ref{prop:dK<Sum}, using the
fact that any random variable can be represented in distribution as a
function of a uniformly distributed random variable, and makes more
precise the central limit theorem of 
   \cite{dejong1987,dejong1990}.
   In comparison with Theorem~1.3 in \cite{doblerpeccati}, 
   Theorem~\ref{thm:S}
   is stated for the Kolmogorov distance instead of the Wasserstein
   distance, it applies to Hoeffding decompositions
   in full generality and not only to degenerate $U$-statistics
   for which $|J|$ is constrained to a fixed value $|J|=d$ for
   some $d\in \{1,\ldots , n\}$ in the sum
   \eqref{eq:Hd}. 
\begin{theorem}\label{thm:S}
  Let $1 \leq d \leq n$. For any $W_n\in L^4(\Omega)$ admitting the Hoeffding decomposition \eqref{eq:Hd} with $|J|\leq d$, and such that $\E\[W_n^2\]=1$,
     we have
     \begin{align}
       \nonumber 
    d_K(W_n,\mathcal N)\leq C_d&\left( \sum_{0\leq l< i\leq d}\sum_{|J|=i-l}\E\Bigg[\Bigg(\sum_{|K|=l, K\cap J=\phi}\E \big[\(W_{J\cup K}\)^2\mid \mathcal F_J\big]\Bigg)^2\Bigg]
    \right.
    \\
       \nonumber 
  &\  +\sum_{1\leq l< i\leq d}
    \sum_{\substack{|J_1|=|J_2|=i-l\\J_1\cap J_2=\phi}}\E\Bigg[\Bigg(\sum_{\substack{|K|=l, K_1\cap\(J_1\cup J_2\)=\phi}}\E\big[W_{J_1\cup K}W_{J_2\cup K} \mid \mathcal F_{J_1\cup J_2}\big]\Bigg)^2\Bigg]
    \\
    \label{bd01}
    &
    \left. +\sum_{1\leq l< i\leq d}\sum_{|J|=i-l}\E\Bigg[\Bigg(\sum_{|K|=l, K\cap J=\phi}\E\big[W_{K}W_{J\cup K} \mid \mathcal F_J\big]\Bigg)^2\Bigg]\right)^{1/2}, 
  \end{align}
  where $C_d > 0$ depends only on $d$.
   \end{theorem} 
    \begin{Proof}
      By representing $X_i$ as
      $X_{i}\stackrel{d}{=}F_i^{-1}\((U_i+1)/2\)$
      where $F_i^{-1}$ is the generalized inverse of the
      cumulative distribution function $F_i$ of $X_i$,
      $i=1,\ldots , n$, we rewrite
\eqref{eq:Hd2} as the sum of multiple stochastic integrals 
      \begin{align*}
        W_n &\stackrel{d}{=}\sum_{k=1}^d I_k(f_k),\end{align*}
where 
\begin{align}\label{eq:fk}
  & f_k(x_1,\ldots ,x_k) :=
  \\\nonumber
  & \frac{1}{k!}
  \sum_{ J = \{i_1,\ldots ,i_k\} \subset [n]} 
  g_J 
  \(
  F_{i_1}^{-1}\(\frac{x_1}{2} -\Big\lfloor \frac{x_1}{2}\Big\rfloor\),
  \ldots ,
   F_{i_k}^{-1}\(\frac{x_k}{2} - \Big\lfloor \frac{x_k}{2} \Big\rfloor \)
  \)
{\bf 1}_{[2i_1-2, 2i_1)\times\cdots \times [2i_k-2,2i_k)}(x_1,\ldots ,x_k),
\end{align}
 $(x_1,\ldots ,x_k)\in \real_+^k$. 
    Next, denoting
    $$
    \widehat{\bf N}^m :=\big\{ ( k_1,\ldots ,k_m) \ : \ k_1,\ldots , k_m \geq 1, \ k_i\neq k_j\text{ if }i\neq j, \ 1\leq i,j\leq m \big\},
    $$
    we have 
  \begin{align*}
& \| f_i\star_{i}^{l}f_i\|^2_{L^2 (\R_+^{i-l})}\\
    &=\frac{1}{2^{2l}}\sum_{\bold j\in \widehat{\bf N}^{i-l}} \int_{[2j_1-2,2j_1)\times\cdots \times[2j_{i-l}-2,2j_{i-l})}
        \left(
        \sum_{\bold k\in \widehat{\bf N}^l}\int_{[2k_1-2,2k_1)\times\cdots \times[2k_{l}-2,2_{l})}(f_i(x_1,\ldots ,x_i))^2dx_1\cdots dx_l\right)^2
            \\
            & \qquad \qquad \qquad \qquad \qquad \qquad \qquad \qquad \qquad \qquad  \qquad \qquad \qquad \qquad \qquad \qquad \qquad \qquad 
            dx_{l+1}\cdots dx_{i}\\
            &\leq (i-l)!(l!)^2\sum_{\substack{|J|=i-l\\J=\{j_1, \ldots ,j_{i-l}\}}}\int_{[2j_1,2j_1+2)\times\cdots \times[2j_{i-l},2j_{i-l}+2)}
                                \\
            & \qquad \qquad\qquad  \qquad \(\sum_{\substack{|K|=l\\K=\{k_1, \ldots ,k_{l}\}}}\int_{[2k_1,2k_1+2)\times\cdots \times[2k_{l},2_{l}+2)}\big(f_i(x_1,\ldots ,x_i)\big)^2dx_1\cdots dx_l\)^2
                                dx_{l+1}\cdots dx_{i}\\
                                &\leq C\sum_{|J|=i-l}\E\Bigg[\Bigg(
                                  \sum_{
                                    |K|=l, K\cap J=\phi}\E\big[\(W_{J\cup K}\)^2\mid \mathcal F_J\big]\Bigg)^2\Bigg],
  \end{align*}
for some $C=C(d)$. Similarly, we get $$\left\| f_{i}\star_{l}^l f_{i}\right\|^2_{
  L^2(\R_+^{2(i-l)})}\leq C\sum_{\substack{|J_1|=|J_2|=l-i\\J_1\cap J_2=\phi}}\E\Bigg[\Bigg(
 \sum_{\substack{|K|=l, K\cap\(J_1\cup J_2\)=\phi}}\E\big[W_{J_1\cup K}W_{J_2\cup K} \mid \mathcal F_{J_1\cup J_2}\big]\Bigg)^2\Bigg]$$
      and 
      $$\left\| f_{l}\star_{l}^l f_{i}\right\|^2_{
        L^2(\R_+^{i-l})}\leq C\sum_{|J|=i-l}\E\Bigg[\Bigg(\sum_{|K|=l, K\cap J=\phi}\E\big[W_{K}W_{J\cup K} \mid \mathcal F_J\big]\Bigg)^2\Bigg],
      \quad 1 \leq l < i \leq d.
      $$
      We conclude by applying the above to Proposition~\ref{prop:dK<Sum},
      which yields the required bound.
    \end{Proof}
    \subsection{Degenerate $U$-statistics}
    In this section we narrow  our attention to the degenerate $U$-statistics of a given order $d\geq1$, which are random variables
    $W_{n,d}$ admitting the Hoeffding decomposition \eqref{eq:Hd} with $|J|=d$.
    \begin{theorem}\label{thm:degenerateU}
     For any degenerate $U$-statistics  $W_{n,d}\in L^4(\Omega)$ of order $d\geq1$, and such that $\E\[W_{n,d}^2\]=1$,
     we have
     \begin{align*}
       \nonumber 
    d_K(W_{n,d},\mathcal N)&\leq\sqrt{\Var\left[\sum_{k=1}^\infty \E\big[ \(W_{n,d}-\E\big[W_{n,d}|\{X_k\}^c\big]\)^2\, | \{X_k\}^c\big]
    \right]}\\
    &\ \ \ +24\sqrt{2\,\E\sum_{k=1}^\infty \E\left[ \big(W_{n,d}-\E\big[W_{n,d}|\{X_k\}^c\big]\big)^4 \right]}\\
    &\leq C_d\Bigg( \sum_{0\leq l< d} \ \sum_{|J|=d-l}\E\Bigg[\Bigg(\sum_{|K|=l, K\cap J=\phi}\E \big[\(W_{J\cup K}\)^2\mid \mathcal F_J\big]\Bigg)^2\Bigg]
    \phantom{\sum_{\substack{|K|=l, K_1\cap\(J_1\cup J_2\)=\phi}}}
    \\
       \nonumber 
  &\  +\sum_{1\leq l<  d}
    \ \sum_{\substack{|J_1|=|J_2|=d-l\\J_1\cap J_2=\phi}}\E\Bigg[\Bigg(\sum_{\substack{|K|=l, K_1\cap\(J_1\cup J_2\)=\phi}}\E\big[W_{J_1\cup K}W_{J_2\cup K} \mid \mathcal F_{J_1\cup J_2}\big]\Bigg)^2\Bigg]\Bigg)^{1/2}, 
  \end{align*}
  where $\{X_k\}^c=\{X_1, \ldots ,X_{k-1},X_{k+1}, \ldots ,X_n \}$ and $C_d > 0$ depends only on $d$.
   \end{theorem} 
    \begin{Proof}
   The first bound is just the latter bound from Proposition \ref{prop:dKI} rewritten in a different form. Namely, it is enough to take $f_d$ as in \eqref{eq:fk} and then we have for $t\in[2k,2k+2)$
   $$\nabla_tW_{n,d}=\E\[W_{n,d}\mid \{X_k\}^c, X_k=t\]-\E\[W_{n,d}\mid \{X_k\}^c\].$$
     The other bound in the assertion follows from Proposition \ref{prop:dKI} in view of \eqref{eq:...<R_X}, \eqref{eq:R_X<} -- where the last sum is vanishing -- and the proof of Theorem \ref{thm:S}.
    \end{Proof}
    \subsubsection*{Weighted $U$-statistics}
    As an example,  we consider classical 
    degenerate weighted $U$-statistics. Precisely,
    given $(X_1,\ldots , X_n)$ an i.i.d. sequence of random variables
    with distribution $\nu$, 
    we define
     \begin{equation}
       \label{und} 
     U_{n,d} ={n\choose d}^{-1}\sum_{1\leq k_1< \cdots < k_d\leq n}w(k_1,\ldots ,k_d)g\(X_{k_1},\ldots ,X_{k_d}\),
     \qquad 1 \leq d \leq n,
\end{equation} 
     where $w(k_1,\ldots ,k_d)\in\R$ is symmetric and
     vanishes on diagonals, and
     $g\(X_{k_1},\ldots ,X_{k_d}\) \in L^2(\Omega)$,
     $1\leq k_1<\cdots < k_d\leq n$,  satisfies 
         \begin{align}
         \label{eq:Ucondition}
    \E\[g\(X_1,x_2,\ldots x_d\)\]=0, \qquad (x_2,\ldots ,x_d)\in\R^{d-1}.
    \end{align}
         The variance $\sigma^2$ of $U_{n,d}$ is given by
          $$\sigma^2:=\Var[U_{n,d}]={n\choose d}^{-2}\|g\|^2_{L^2(\R^d,\nu^{\otimes d})}\sum_{1\leq k_1<\cdots <k_d\leq n}w^2(k_1,\ldots ,k_d). 
$$
    The assumption \eqref{eq:Ucondition}  plays a technical role,
    which helps in simplifying the derivations. 
    Nevertheless, it covers important examples of $U$-statistics
    such as quadratic forms and their multidimensional generalizations.
    Sharp bounds have been provided in \cite{CS}
    in case \eqref{eq:Ucondition} is not satisfied, but
    only in the case of classical (i.e. non-weighted) $U$-statistics.
    See also \cite{reichenbachs} for weighted first order $U$-statistics
    based on symmetric Rademacher sequences, and 
    \cite{poly} for a fourth moment type central limit theorem
    in case $g(x_1,\ldots , x_n)=x_1\cdots x_n$
    and $X_1$ has a vanishing third moment. 
\begin{theorem}\label{thm:Ustat}
    Let $U_{n,d}$ be a degenerate weighted $U$-statistics of the form 
    \eqref{und}.
     We have 
    \begin{align*}
    d_K\(\frac{U_{n,d}}{\sigma},\mathcal N \)\leq C_d
    \frac{\|g\|^2_{L^4(\R^d,\nu^{\otimes d})}}{\|g\|^2_{L^2(\R^d,\nu^{\otimes d})}}
    {\frac{\sup_{1\leq l\leq d-1} \sqrt{ \sum_{\bold k, \bold r\in \N^{d-l}}\(\sum_{\bold m\in \N^l}w(\bold k, \bold m)w(\bold r, \bold m)\)^2}}{\sum_{\bold m\in \N^{d}}w^2(\bold m)}}
    \end{align*} 
    for some $C_d>0$ depending only on $d=1,\ldots , n$,
    where $\nu$ denotes the distribution of $X_1$. 
    \end{theorem}
        \begin{Proof}
          By Theorem \ref{thm:degenerateU},
                  we have
     \begin{align*}
         &d_K\(\frac{U_{n,d}}{\sigma},\mathcal N\)\\
         &\leq
       \frac{C_d}{\sigma^2} {{n\choose d}^{-2}}
       \Bigg(  \sum_{0\leq l\leq d-1}\int_{\R^{d-l}}\(\int_{\R^l}g^2(x,y) \nu^{\otimes l}(dx)\)^2 \nu^{\otimes(d-l)}(dy)\sum_{\bold k\in\N^{d-l}}\(\sum_{\bold m\in \N^l}w^2(\bold k, \bold m)\)^2
       \\
       &\hspace{25mm} +\sum_{1\leq l\leq d-1}
       \int_{\R^{d-l}} \int_{\R^{d-l}}\(\int_{\R^l}g(x,y)g(x,z) \nu^{\otimes l}(dx)\)^2 \nu^{\otimes(d-l)}(dy) \nu^{\otimes(d-l)}(dz)
       \\
       &\hspace{70mm}\times\sum_{\bold k, \bold r\in\N^{d-l}}\(\sum_{\bold m\in \N^l}w(\bold k, \bold m)w(\bold r, \bold m)\)^2
     \Bigg)^{1/2}. 
  \end{align*}
  Since $\nu$ is a probability measure, we have
  \begin{eqnarray*}
     \int_{\R^{d-l}}\(\int_{\R^l}g^2(x,y) \nu^{\otimes l}(dx)\)^2 \nu^{\otimes(d-l)}(dy)
     & \leq & \int_{\R^{d-l}}\int_{\R^l}g^4(x,y) \nu^{\otimes l}(dx) \nu^{\otimes(d-l)}(dy)
     \\
      & = & \|g\|^4_{L^4(\R^d,\nu^{\otimes d})},
\end{eqnarray*} 
  as well as
  \begin{align*}
    &
    \! \! \! \! \! \! \! \! \! \! \! \! \! \! \! \! \! \! 
    \int_{\R^{d-l}}\int_{\R^{d-l}}\(\int_{\R^l}g(x,y)g(x,z) \nu^{\otimes l}(dx)\)^2 \nu^{\otimes(d-l)}(dy) \nu^{\otimes(d-l)}(dz)\\
    &\leq \int_{\R^{d-l}}\int_{\R^{d-l}}\int_{\R^l}g^2(x,y) \nu^{\otimes l}(dx)\int_{\R^l}g^2(x,z) \nu^{\otimes l}(dx) \nu^{\otimes(d-l)}(dy) \nu^{\otimes(d-l)}(dz)
    \\
    &=\|g\|^4_{L^2(\R^d,\nu^{\otimes d})}\leq \|g\|^4_{L^4(\R^d,\nu^{\otimes d})}.
  \end{align*}
  Using also the inequality
  \begin{align*}
  \sum_{\bold k\in\N^{d-l}}\(\sum_{\bold m\in \N^l}w^2(\bold k, \bold m)\)^2\leq \sum_{\bold k, \bold r\in\N^{d-l}}\(\sum_{\bold m\in \N^{l}}w( \bold k, \bold m)w( \bold r, \bold m)\)^2,
  \end{align*}
we arrive at
 \begin{align*}
 d_K\(\frac{U_{n,d}}{\sigma},\mathcal N\)&\leq \frac{C_d}{\sigma^2{{n\choose d}^2}}\|g\|^2_{L^4(\R^d,\nu^{\otimes d})}\sup_{1\leq l\leq d-1}\sqrt{  \sum_{\bold k, \bold r\in\N^{d-l}}\(\sum_{\bold m\in \N^{l}}w( \bold k, \bold m)w( \bold r, \bold m)\)^2}\\[10pt]
 &= {C_d}\frac{\|g\|^2_{L^4(\R^d,\nu^{\otimes d})}}{\|g\|^2_{L^2(\R^d,\nu^{\otimes d})}}\frac{\sup_{1\leq l\leq d-1} \sqrt{ \sum_{\bold k, \bold r\in\N^{d-l}}\(\sum_{\bold m\in \N^l}w(\bold k, \bold m)w(\bold r, \bold m)\)^2}}{\sum_{1\leq k_1, \ldots ,k_d\leq n}w^2(k_1,\ldots ,k_d)},
 \end{align*}
 which is  the  bound in the assertion. 
      \end{Proof}
      
\subsection{Random graphs}
Consider the 
 \cite{ER}
    random graph $\mathbb{G}_n(p)$ 
constructed by independently retaining any edge in 
the complete graph $K_n$ on $n$ vertices with probability $p\in(0,1)$.
Here, we assign 
an independent sample of a random weight $X$
to every edge in $\mathbb{G}_n (p_n)$,
and we consider the renormalized random weight 
$$\widetilde{W}^G_n:=\frac{W^G_n-\E[W^G_n]}{\sqrt{\Var[W^G_n]}},
$$
where $W^G_n$ denotes the combined weight of graphs
in $\mathbb{G}_n(p_n)$ that are isomorphic to
a fixed graph $G$.  
By writing the combined weight $W^G_n$ of graphs
in $\mathbb{G}_n(p_n)$ that are isomorphic to
a fixed graph $G$  as a sum of multiple stochastic integrals (which is equivalent to finding its Hoeffding decomposition)
we obtain the following result as in
\cite{PS3}, by replacing the use of Thereom~5.1 therein
with Theorem~\ref{prop:dK<Sum} above. 
\begin{theorem}
      \label{rg}
      Let  $G$ be a graph without isolated vertices.
 The renormalized weight $\widetilde{W}^G_n$
 of graphs in $\mathbb{G}_n(p_n)$ that are isomorphic to
 $G$ satisfies 
\begin{equation} 
\nonumber 
  d_K \big(\widetilde{W}^G_n,\mathcal{N}\big) \leq C
  \frac{\sqrt{\E\[\(X-\E[X]\)^4\]}+(1-p)(\E[X])^2}{\Var [X]+(1-p)(\E[X])^2}
  \((1-p)\min_{\substack{ H\subset G\\e_H\geq1}} n^{v_H}p^{e_H} \)^{-1/2},
\end{equation}
for some constant  $C=C(e_G)>0$.
\end{theorem}
    This extends other
    Kolmogorov distance bounds previously obtained for triangle counting
in \cite{nross}, and in \cite{reichenbachsAoP}
 using the Malliavin approach to 
 the Stein method, see also \cite{roellin2} for triangle counting and \cite{PS2}
 for arbitrary subgraph counting, 
 and \cite{reichenbachs} for weighted first order Rademacher $U$-statistics
 in the symmetric case $p=1/2$. 
As a consequence, if $p_n$ satisfies $p_n<c<1$, $n\geq 1$, we have
\begin{align*}
d_K \big(\widetilde{W}^G_n,\mathcal{N}\big)&\leq C\frac{\sqrt{\E\[X^4\]}}{\E\[X^2\]}\((1-p_n)
  \min_{\substack{ H\subset G\\e_H\geq1}} n^{v_H}p_n^{e_H} \)^{-1/2}, 
\end{align*}
 and for  $p_n>c>0$, $n\geq 1$, it holds
\begin{equation}
\nonumber 
d_K \big(\widetilde{W}^G_n,\mathcal{N}\big)\leq C
\frac{\sqrt{\E\[X^4\]}}{
  n\sqrt{1-p_n} \Var [X] }.
 \end{equation} 
 Applications to cycle graphs,
 complete graphs
 trees can be treated as in \cite{PS3} by replacing the
 Kolmogorov distance with the Wasserstein distance.
 
\section{Quadratic forms}
 \label{section:qf}
 We consider the quadratic form $Q_n$ defined as  
 $$
 Q_n=\sum_{\substack{1\leq i,j\leq n\\i\neq j}}a_{ij}X_i X_j
 + \sum_{k=1}^n a_{kk} \(X_k^2-\E\[X_k^2\]\),
 $$
 where $A_n= ( a_{ij} )_{1\leq i,j\leq n}$
    is a symmetric matrix, $n\geq 1$,
    and $\(X_k\)_{k\geq1}$ denotes
    i.i.d. copies of a given random variable $X$ satisfying $\E\[X\]=0$.
    In the sequel, we let $\mu_k:=\E[X^k]$,
 $\tilde{\mu}_k:=\E\big[\(X^2-\E\[X^2\]\)^{k/2}\big]$,
 $k\geq 2$,
    and
    $$
\sigma_n^2
  := 
  \Var [ Q_n] = \E\[Q_n^2\]
  =2 \mu_2^2 
  \sum_{\substack{1\leq i,j\leq n\\i\neq j}}a_{ij}^2 
  + \tilde\mu_4 \sum_{i=1}^na_{ii}^2.
  $$ 
  Many papers in the literature are devoted to asymptotical normality of quadratic forms.
  The best known convergence rates in the general case where the diagonal
  of $A$ may not vanish are given in 
     \cite{goetze-tikhomirov2002}, as 
\begin{align}\label{eq:GT}
 d_K\(\frac{Q_n}{\sigma_n},\mathcal{N}\)\leq C(\gamma)\frac{\(\E\[|X|^3\]\)^2+\gamma \E[X^6]}{\sqrt{\sum_{1\leq i,j\leq n}a_{ij}^2}}|\lambda_1|,
\end{align}
see Theorem~1.1 therein, 
where $\lambda_1$ denotes the largest 
     absolute eigenvalue of $A_n$, $\gamma=\sum_{i=1}^n a_{ii}^2/\sum_{1\leq i,j\leq n}a_{ij}^2$,
     and the constant $C(\gamma)$ blows up when $\gamma$ tends to one, i.e. when the linear part is dominating.
     \subsubsection*{Vanishing diagonals}
     More is known if we assume the diagonal of $A_n$ to be empty,
     in which case 
     \cite{dejong1987} proved the asymptotic normality of $Q_n/\sigma_n$ under the conditions  
    \begin{align}\label{eq:dJ}\E\[\({Q_n}/{\sigma_n}\)^4\]\longrightarrow 3\ \text{ ~and~ } \frac{1}{\sigma_n^2} \max_{1\leq i \leq n}\sum_{j=1}^na_{ij}^2\longrightarrow0.\end{align}
 In addition, for $(X_k)_{k\geq1}$ a Rademacher sequence, Theorem~1.1 
 in \cite{dobler} restricted to double integrals gives the corresponding bound 
\begin{align}\label{eq:DK}
d_K\(\frac{Q_n}{\sigma_n},\mathcal{N}\)\leq C\(\sqrt{\big|\E\[(Q_n/\sigma_n )^4\]-3\big|}+ \frac{1}{\sigma_n} \sqrt{\max_{1\leq i \leq n}\sum_{1\leq j\leq 1}a_{ij}^2}\).
\end{align}
The same bound may be concluded from \cite{doblerpeccati} for $(X_k)_{k\geq1}$ being any i.i.d. sequence, but only in Wasserstein distance. Note that the quantity $\displaystyle \max_{1\leq i \leq n}\sum_{j=1}^na_{ij}^2$ corresponds to ``maximal influence'', see
 \cite{mossel}, \cite{nourdin4}.

 \medskip

 The bound 
 \begin{equation}
   \label{fjkdslfd} 
      d_W\(\frac{Q_n}{\sigma_n},\mathcal N \)\leq
      C \frac{\mu_4}{\sigma_n^2}\(
      \sqrt{\sum_{i=1}^n\Bigg(\sum_{k=1}^n a_{ik}^2\Bigg)^2}+
      \sqrt{\sum_{i, j=1}^n\Bigg(\sum_{k=1}^n a_{ik}a_{kj}\Bigg)^2}\).
\end{equation} 
 has been provided for Rademacher sequences 
 using the Wasserstein distance in Proposition~3.1 of \cite{chatterjee},
 and has been recently extended
 to arbitrary i.i.d. sequences using the Kolmogorov distance in 
 \cite{shao2}, Theorem~3.1.

 \medskip
 
 Corollary~\ref{chkskl} recovers this bound as
 an immediate consequence of Theorem~\ref{thm:Ustat} by taking $d=2$,
    $w(k_1,k_2)=a_{k_1k_2}$, $1\leq k_1,k_2\leq n$, $k_1\not= k_2$, 
    and $g(y_1,y_2)=y_1 y_2$.
 Note however that 
 only the second term
 is significant  in the right-hand side of \eqref{fjkdslfd},
 making the conjecture at the end of Section~3.1 in  \cite{shao2}
 pointless.
\begin{corollary}
  \label{chkskl}
  Assume  $a_{ii}=0$, $i=1,\ldots , n$.
  Then, there exists a constant $C>0$ such that
    \begin{align*}
      d_K\(\frac{Q_n}{\sigma_n},\mathcal N \)\leq
    C \frac{\mu_4}{\sigma_n^2}
 \sqrt{ \sum_{i, j=1}^n\Bigg(\sum_{k=1}^n a_{ik}a_{kj}\Bigg)^2}= C \frac{\mu_4}{\sigma_n^2}
 \sqrt{ {\rm Tr}(A_n^4) }, \qquad n \geq 1.  
    \end{align*}
\end{corollary}
Corollary~\ref{chkskl} also improves \eqref{eq:GT} for matrices $A_n$ with empty diagonal, since 
\begin{align}
  \label{eq:Tr<l}
    \sqrt{{\rm Tr}(A_n^4)}
    =
    \sqrt{\sum_{k=1}^n \lambda_k^4}
    \leq
    |\lambda_1| \sqrt{\sum_{k=1}^n \lambda_k^2}
    \leq 
    |\lambda_1| \sqrt{\sum_{i,j=1}^n a_{ij}^2}
    \leq 
    \frac{\sigma_n}{\mu_2} |\lambda_1|.
   \end{align}
\subsubsection*{Non-empty diagonals} 
Theorem~\ref{qf} below generalizes and improves all the aforementioned results.
 First, in comparison with the above bound \eqref{eq:GT} 
 of \cite{goetze, goetze-tikhomirov2002}, 
 it gives better rates under weaker assumptions, as noted in \eqref{eq:Tr<l}. 
 Furthermore, it extends every other result
 by applying as well to non-vanishing diagonals.
 In addition, it completes Corollary \ref{chkskl} with an
 additional bound related to so called fourth moment phenomenon
 (\cite{nualart2005central}), 
 and it also extends \eqref{eq:DK}
 from the Rademacher case to any distribution. 
 Finally, it deals with the Kolmogorov distance instead
 of the Wasserstein distance considered in \cite{doblerpeccati}. See also Theorem 3.11 in \cite{BC} for some bounds in total variation and Kolmogorov distances, which however provide worse rates and require slightly stronger assumptions.
\begin{theorem}
\label{qf}
 There exist absolute constants $C_1,C_2>0$ such that
\begin{equation} 
  \label{r1}
   d_K\(\frac{Q_n}{\sigma_n} ,\mathcal N \)
 \leq {C_1} \(\sqrt{\big|\E\[(Q_n/\sigma_n )^4\]-3\big|}+ \frac{\alpha_n} {\sigma_n}\sqrt{\max_{1\leq i \leq n}\sum_{1\leq j\leq 1}a_{ij}^2}\),
\end{equation} 
and 
\begin{equation} 
  \label{r2}
  d_K\(\frac{Q_n}{\sigma_n} ,\mathcal N \)\leq C_2 \frac{\beta_n}{\sigma_n^2}
 \sqrt{ {\rm Tr}(A_n^4) },
\end{equation} 
where
$$
\alpha_n :=\mu_2+\frac{\mu_4}{\mu_2}\mathbbm{1}_{\big\{ a_{11}^2+\cdots + a_{nn}^2 >0 \big\}}, \ \ \ \text{ and }\ \ \  \beta_n=\mu_4+\sqrt{\mu_8}\,\mathbbm{1}_{\big\{ a_{11}^2+\cdots + a_{nn}^2 >0 \big\}}.$$
\end{theorem}

\begin{Proof}
 The quadratic form $Q_n$ admits the Hoeffding decomposition
$$
Q_n=\sum_{1\leq i,j \leq n}W_{\{i,j\}}+\sum_{k=1}^n W_{\{k\}},$$
where 
\begin{align*}
  W_{\{i,j\}}&=2a_{ij}X_iX_j,
  \qquad
  W_{\{k\}}=a_{kk}\(X_k^2-\E\[X_k^2\]\).
\end{align*}
Thus, Theorem \ref{thm:S} gives 
\begin{align} \label{eq:bound1}
  d_K\(\frac{Q_n}{\sigma_n}, \mathcal N\)\leq \frac{C}{\sigma_n^2}&\Bigg(
  \tilde{\mu}_8\sum_{i=1}^na_{ii}^4+2\mu_4^2\sum_{\substack{1\leq i,j\leq n\\i\neq j}}a_{ij}^4+2\mu_2^2\mu_4\sum_{\substack{1\leq i, j, k\leq n\\ i\neq j, i\neq k, j\neq k}}a_{ij}^2a_{ik}^2
  \\
  \nonumber 
& \ \  +\mu_2^4\sum_{\substack{1\leq i,j\leq n\\i\neq j}}\Bigg(\sum_{\substack{1\leq k\leq n\\k\neq i,j}}a_{ik}a_{kj}\Bigg)^2+\mu_3^2\mu_2\sum_{i=1}^n\Bigg(\sum_{\substack{1\leq j\leq n\\i\neq j}}a_{jj}a_{ij}\Bigg)^2\Bigg)^{1/2}.
\end{align}
 Next, we estimate this bound by means of $\E\[Q_n^4\]$ and
$\displaystyle \max_{1\leq i\leq n}\sum_{\substack{1\leq j\leq n}}a_{ij}^2$.
A direct calculation shows that
$$\E\[Q_n^4\]=S_1+3S_2+4S_3,$$
where
\begin{align*}
S_1&:=\tilde{\mu}_8\sum_{i=1}^na_{ii}^4+16\mu_4^2\sum_{1\leq i<j\leq n}a_{ij}^4+48\mu_2^2\mu_4\sum_{\substack{1\leq i, j, k\leq n\\ i\neq j, i\neq k, j\neq k}}a_{ij}^2a_{ik}^2\\
&\ \ \ +48\mu_2^4\sum_{\substack{1\leq i_1,i_2,i_3,i_4\leq n\\i_k\neq i_l \text{ if }k\neq l}}a_{i_1i_2}a_{i_2i_3}a_{i_3i_4}a_{i_4i_1}+48\mu_3^2\mu_2\sum_{\substack{1\leq i, j, k\leq n\\ i\neq j, i\neq k, j\neq k}}a_{ii}a_{jj}a_{ik}a_{kj}\\
&\ \ \ +48\mu_3^2\mu_2\sum_{\substack{1\leq i, j, k\leq n\\ i\neq j, i\neq k, j\neq k}}a_{kj}^2a_{ik}a_{ij} ,
\end{align*} 
and
 $$ 
S_2:=\tilde{\mu}_4^2\sum_{i\neq j}a^2_{ii}a^2_{jj}+4\tilde{\mu}_4\mu_2^2\sum_{\substack{1\leq i,j,k\leq n\\j\neq k,\ j,k\neq i }}a_{ii}^2a^2_{jk}+ 4\mu_2^4\sum_{\substack{1\leq i_1,i_2,i_3,i_4\leq n\\i_k\neq i_l \text{ if }k\neq l}}a_{i_1i_2}^2a^2_{i_3i_4},
$$
and
\begin{align*}
S_3&:=3\tilde{\mu}_4^2\sum_{\substack{1\leq i,j\leq n\\i\neq j}}a_{ii}a_{jj}a_{ij}^2+8\mu_3\(\mu_5-\mu_3\mu_2\)\sum_{i\neq j}a_{ii}a_{ij}^3+6\mu_2\(\tilde{\mu}_6+\tilde{\mu}_4\mu_2\)\sum_{\substack{1\leq i,j\leq n\\i\neq j }}a_{ii}^2a^2_{ij}\\
&\ \ \ \ +12\mu_3^2\mu_2\sum_{\substack{1\leq i,j,k\leq n\\i\neq j, j\neq k, i\neq k}}a_{ii}a_{ij}a^2_{jk}+24\mu_2^2\tilde{\mu}_4\sum_{\substack{1\leq i,j,k\leq n\\i\neq j, j\neq k, i\neq k}}a_{ii}a_{ij}a_{ik}a_{kj}.
\end{align*}
 The sum $S_1$ is supposed to dominate the right-hand side of \eqref{eq:bound1}, $S_2$ is approximating $\sigma^2$, and $S_3$  contains some remainders that are problematic due to their unknown sign and vanishes if the diagonal of $A$ is empty.  First, by 
$$ \sum_{\substack{1\leq i,j\leq n\\i\neq j}}\Bigg(\sum_{\substack{1\leq k\leq n\\k\neq i,j}}a_{ik}a_{kj}\Bigg)^2=  \sum_{\substack{1\leq i_1,i_2,i_3,i_4\leq n\\i_k\neq i_l \text{ if }k\neq l}}a_{i_1i_2}a_{i_2i_3}a_{i_3i_4}a_{i_4i_1}+\sum_{\substack{1\leq i, j, k\leq n\\ i\neq j, i\neq k, j\neq k}}a_{ik}^2a_{kj}^2$$
and
$$\sum_{i=1}^n\Bigg(\sum_{\substack{1\leq j\leq n\\i\neq j}}a_{jj}a_{ij}\Bigg)^2
= \sum_{\substack{1\leq i, j, k\leq n\\ i\neq j, i\neq k, j\neq k}}a_{ii}a_{jj}a_{ik}a_{kj}+\sum_{\substack{1\leq i,j\leq n\\i\neq j }}a_{ii}^2a^2_{ij},$$
we get
\begin{align*}
S_1&:=\tilde{\mu}_8\sum_{i=1}^na_{ii}^4+16\mu_4^2\sum_{1\leq i<j\leq n}a_{ij}^4+48\mu_2^2\mu_4\sum_{\substack{1\leq i, j, k\leq n\\ i\neq j, i\neq k, j\neq k}}a_{ij}^2a_{ik}^2\\
&\ \ \ +48\mu_2^4\sum_{\substack{1\leq i,j\leq n\\i\neq j}}\Bigg(\sum_{\substack{1\leq k\leq n\\k\neq i,j}}a_{ik}a_{kj}\Bigg)^2+48\mu_3^2\mu_2\sum_{i=1}^n\Bigg(\sum_{\substack{1\leq j\leq n\\i\neq j}}a_{jj}a_{ij}\Bigg)^2\\
&\ \ \ -48\mu_2^4\sum_{\substack{1\leq i, j, k\leq n\\ i\neq j, i\neq k, j\neq k}}a_{ik}^2a_{kj}^2-48\mu_3^2\mu_2 \sum_{\substack{1\leq i,j\leq n\\i\neq j }}a_{ii}^2a^2_{ij}+48\mu_3^2\mu_2\sum_{\substack{1\leq i, j, k\leq n\\ i\neq j, i\neq k, j\neq k}}a_{kj}^2a_{ik}a_{ij}.
\end{align*} 
The first two lines dominate  the right-hand side of \eqref{eq:bound1} with substantial surplus, which will be used to deal with the last term of $S_1$ and some  terms of $S_3$. Indeed, by $\mu_3^2\mu_2\leq \mu_4\mu_2^2$ and the inequality of arithmetic and geometric means, we have
\begin{align}\nonumber
& 48\mu_3^2 \mu_2 \left|\sum_{\substack{1\leq i, j, k\leq n\\ i\neq j, i\neq k, j\neq k}}a_{kj}^2a_{ik}a_{ij}\right|\\\nonumber
  &\leq 46 \mu_4\mu_2^2\sum_{\substack{1\leq i, j, k\leq n\\ i\neq j, i\neq k, j\neq k}}\frac12\Big((a_{kj}a_{ik})^2+(a_{kj}a_{ij})^2\Big)
  +\sum_{\substack{1\leq k, j\leq n\\ k\neq j}} \Bigg((\mu_4a_{kj}^2)^2+\Bigg(\mu_2^2\sum_{\substack{1\leq i\leq n\\ i\neq j,k}}a_{ik}a_{ij}\Bigg)^2\Bigg)\\\label{eq:S1a}
&=46\mu_2^2\mu_4\sum_{\substack{1\leq i, j, k\leq n\\ i\neq j, i\neq k, j\neq k}}a_{ij}^2a_{ik}^2+\mu_4^2\sum_{1\leq i<j\leq n}a_{ij}^4+\mu_2^4\sum_{\substack{1\leq i,j\leq n\\i\neq j}}\Bigg(\sum_{\substack{1\leq k\leq n\\k\neq i,j}}a_{ik}a_{kj}\Bigg)^2.
\end{align}
Since, additionally
\begin{align}\nonumber
  & \left|\mu_2^4\sum_{\substack{1\leq i, j, k\leq n\\ i\neq j, i\neq k, j\neq k}}a_{ik}^2a_{kj}^2+\mu_3^2\mu_2 \sum_{\substack{1\leq i,j\leq n\\i\neq j }}a_{ii}^2a^2_{ij}\right|
  \\
  \label{eq:S1b}
&\leq \sigma_n^2 \mu_2^2\max_{1\leq i \leq n}\sum_{j=1}^na_{ij}^2+\sigma_n^2\frac{\mu_3^2}{\mu_2}
\mathbbm{1}_{\big\{ a_{11}^2+\cdots + a_{nn}^2 > 0\big\}}
\max_{1\leq i \leq n}\sum_{j=1}^na_{ij}^2
 \leq \sigma_n^2\alpha_n^2\max_{1\leq i \leq n}\sum_{j=1}^na_{ij}^2,
\end{align}
 we arrive at
 \begin{eqnarray}
   \nonumber
   \lefteqn{
     \! \! \! \! \! \! \! \! \! \! 
     d_K\(\frac{Q_n}{\sigma_n}, \mathcal N\) 
  \leq  \frac{C}{\sigma_n^2}\Bigg(
  S_1+48\sigma_n^2\alpha_n^2\max_{1\leq i \leq n}\sum_{j=1}^na_{ij}^2
   }
   \\
  \nonumber
& & +46\mu_2^4\sum_{\substack{1\leq i,j\leq n\\i\neq j}}\Bigg(\sum_{\substack{1\leq k\leq n\\k\neq i,j}}a_{ik}a_{kj}\Bigg)^2+47\mu_3^2\mu_2\sum_{i=1}^n\Bigg(\sum_{\substack{1\leq j\leq n\\i\neq j}}a_{jj}a_{ij}\Bigg)^2\Bigg)^{1/2}\\\nonumber
  &\leq & \frac{C}{\sigma_n^2} \Bigg( \E\[Q_n^4\]-3\sigma_n^4+ 48 \sigma_n^2 \alpha_n^2\max_{1\leq i \leq n}\sum_{j=1}^na_{ij}^2+3 (\sigma_n^4-S_2)
  \\
  \label{eq:bound2}
  & &  -4 S_3
  -24\mu_2^4\sum_{\substack{1\leq i,j\leq n\\i\neq j}}\Bigg(\sum_{\substack{1\leq k\leq n\\k\neq i,j}}a_{ik}a_{kj}\Bigg)^2
  -24\mu_3^2\mu_2\sum_{i=1}^n\Bigg(\sum_{\substack{1\leq j\leq n\\i\neq j}}a_{jj}a_{ij}\Bigg)^2
  \Bigg)^{1/2}.\quad 
\end{eqnarray} 
 Next, in order to bound $3(\sigma_n^4-S_2)$, we calculate
\begin{align*}
\sigma_n^4&=\Bigg(2\mu_2^2\sum_{\substack{1\leq i,j\leq n\\i\neq j}}a_{ij}^2+\tilde{\mu}_4\sum_{i=1}^na_{ii}^2\Bigg)^2\\
&=4\mu_2^4\Bigg(\sum_{\substack{1\leq i,j\leq n\\i\neq j}}a_{ij}^2\Bigg)^2+\tilde{\mu}_4^2\Bigg(\sum_{i=1}^na_{ii}^2\Bigg)^2+4\mu_2^2\tilde{\mu}_4\Bigg(\sum_{\substack{1\leq i,j\leq n\\i\neq j}}a_{ij}^2\Bigg)\(\sum_{i=1}^na_{ii}^2\)\\
&=  4\mu_2^4\Bigg(2 \sum_{\substack{1\leq i,j\leq n\\i\neq j}}a_{ij}^4
+\sum_{\substack{1\leq i_1,i_2,i_3,i_4\leq n\\i_k\neq i_l \text{ if }k\neq l}}a_{i_1i_2}^2a^2_{i_3i_4}+2\sum_{\substack{1\leq i,j,k\leq n\\j\neq k,\ j,k\neq i }}a^2_{ij}a^2_{ik}\Bigg)\\
&\ \ \ +\tilde{\mu}_4^2\Bigg(\sum_{i=1}^na_{ii}^4+\sum_{\substack{1\leq i,j\leq n\\i\neq j}}^na_{ii}^2a_{jj}^2\Bigg)+4\mu_2^2\tilde{\mu}_4\Bigg(2\sum_{\substack{1\leq i,j\leq n\\i\neq j}}a_{ii}^2a_{ij}^2+\sum_{\substack{1\leq i,j,k\leq n\\j\neq k,\ j,k\neq i }}a_{ii}^2a^2_{jk}\Bigg),
\end{align*}
hence
\begin{align}\nonumber
  3 \left|\sigma_n^4-S_2\right|&= 24 \mu_2^4\Bigg(\sum_{\substack{1\leq i,j\leq n\\i\neq j}}a_{ij}^4+\sum_{\substack{1\leq i,j,k\leq n\\j\neq k,\ j,k\neq i }}a^2_{ij}a^2_{ik}\Bigg)
  + 3 \tilde{\mu}_4^2\sum_{i=1}^na_{ii}^4+ 24 \mu_2^2\tilde{\mu}_4\sum_{\substack{1\leq i,j\leq n\\i\neq j}}a_{ii}^2a_{ij}^2\\\nonumber
&\leq \max_{1\leq i\leq n}\sum_{\substack{1\leq j\leq n}}a_{ij}^2\Bigg(48\mu_2^4\sum_{\substack{1\leq i,j\leq n\\i\neq j}}a_{ij}^2+27\mu_4\tilde\mu_4\mathbbm{1}_{\big\{ a_{11}^2+\cdots + a_{nn}^2 > 0\big\}}\sum_{i=1}^na_{ii}^2\Bigg)\\\label{eq:S_2}
&\leq \sigma_n^2\(48\mu_2^2+27\frac{\mu_4^2}{\mu_2^2} \mathbbm{1}_{\big\{ a_{11}^2+\cdots + a_{nn}^2 > 0\big\}}\)\max_{1\leq i\leq n}\sum_{\substack{1\leq j\leq n}}a_{ij}^2.
\end{align}
 Regarding $S_3$, we have
\begin{align*}
  \left| \tilde\mu_4^2
  \sum_{\substack{1\leq i,j\leq n\\i\neq j}}a_{ii}a_{jj}a_{ij}^2\right|\leq
  \tilde\mu_4^2\sum_{\substack{1\leq i,j\leq n\\i\neq j}}a_{ii}^2a_{ij}^2
  \leq \sigma_n^2
  \max_{1\leq i\leq n}\sum_{\substack{1\leq j\leq n}}a_{ij}^2,
\end{align*}
and 
\begin{align*}
  & 8 \mu_3\(\mu_5-\mu_3\mu_2\)\sum_{\substack{1\leq i,j\leq n\\i\neq j }} a_{ii}a_{ij}^3
  +6\mu_2\(\tilde{\mu}_6+\tilde{\mu}_4\mu_2\)\sum_{\substack{1\leq i,j\leq n\\i\neq j }}a_{ii}^2a^2_{ij}\\
&= 2 \E\[4\sum_{\substack{1\leq i,j\leq n\\i\neq j }} a_{ii}a_{ij}^3\big(X_i^2-\E\[X_i^2\]\big)X_i^3X_j^3+3\sum_{\substack{1\leq i,j\leq n\\i\neq j }}a_{ii}^2a^2_{ij}\big(X_i^2-\E\[X_i^2\]\big)^2X_i^2X_j^2\]\\
&= 6 \E\[\sum_{\substack{1\leq i,j\leq n\\i\neq j }}\( a_{ii}\big(X_i^2-\E\[X_i^2\]\big)+\frac{2}{3}  a_{ij}X_iX_j\)^2a^2_{ij}X_i^2X_j^2\]-\frac{8}3\E\[\sum_{\substack{1\leq i,j\leq n\\i\neq j }}a^4_{ij}X_i^4X_j^4\]\\
  &\geq -\frac{8}3\mu_4^2\sum_{\substack{1\leq i,j\leq n\\i\neq j }} a^4_{ij}\geq
  -\frac{4}{3} \sigma_n^2\(\frac{\mu_4}{\mu_2}\)^2\max_{1\leq i\leq n}\sum_{\substack{1\leq j\leq n}}a_{ij}^2.
\end{align*}
Furthermore, using the correction terms from \eqref{eq:bound2}, we get
\begin{align*}
  & 12 \mu_3^2\mu_2\sum_{\substack{1\leq i,j,k\leq n\\i\neq j, j\neq k, i\neq k}}a_{ii}a_{ij}a^2_{jk}
  + 6 \mu_3^2\mu_2\sum_{i=1}^n\Bigg(\sum_{\substack{1\leq j\leq n\\i\neq j}}a_{jj}a_{ij}\Bigg)^2\\
&= 6 \mu_3^2\mu_2\sum_{i=1}^n\Bigg(\sum_{\substack{1\leq k\leq n\\k\neq i}}a_{ik}^2+\sum_{\substack{1\leq j\leq n\\j\neq i}}a_{jj}a_{ij}\Bigg)^2- 6 \mu_3^2\mu_2\sum_{i=1}^n\Bigg(\sum_{\substack{1\leq k\leq n\\k\neq i}}a_{ik}^2\Bigg)^2\\
&\geq - 6 \mu_3^2\mu_2\(\max_{1\leq i\leq n}\sum_{\substack{1\leq j\leq n}}a_{ij}^2\)\sum_{i=1}^n\sum_{\substack{1\leq k\leq n\\k\neq i}}a_{ik}^2\geq - 3 \sigma_n^2\(\frac{\mu_4}{\mu_2}\)^2\max_{1\leq i\leq n}\sum_{\substack{1\leq j\leq n}}a_{ij}^2,
\end{align*}
as well as
\begin{align*}
& 24 \mu_2^2\tilde{\mu}_4\sum_{\substack{1\leq i,j,k\leq n\\i\neq j, j\neq k, i\neq k}}a_{ii}a_{ij}a_{ik}a_{kj} + 6 \mu_2^4\sum_{\substack{1\leq i,j\leq n\\i\neq j}}\Bigg(\sum_{\substack{1\leq k\leq n\\k\neq i,j}}a_{ik}a_{kj}\Bigg)^2\\
&= 6 \sum_{\substack{1\leq i,j\leq n\\i\neq j}}\Bigg(2\tilde\mu_4a_{ii}a_{ij}+\mu_2^2\sum_{\substack{1\leq k\leq n\\k\neq i,j}}a_{ik}a_{kj}\Bigg)^2-24\tilde\mu_4^2\sum_{\substack{1\leq i,j\leq n\\i\neq j}}a_{ii}^2a_{ij}^2\\
&\geq - 24 \tilde\mu_4^2\(\max_{1\leq i\leq n}\sum_{\substack{1\leq j\leq n}}a_{ij}^2\) \sum_{\substack{1\leq i,j\leq n\\i\neq j}}a_{ii}^2\geq - 24 \sigma_n^2\mu_4\max_{1\leq i\leq n}\sum_{\substack{1\leq j\leq n}}a_{ij}^2.
\end{align*}
Hence, we arrive at
\begin{align*}
&S_3+6\mu_3^2\mu_2\sum_{i=1}^n\Bigg(\sum_{\substack{1\leq j\leq n\\i\neq j}}a_{jj}a_{ij}\Bigg)^2+ 6 \mu_2^4\sum_{\substack{1\leq i,j\leq n\\i\neq j}}\Bigg(\sum_{\substack{1\leq k\leq n\\k\neq i,j}}a_{ik}a_{kj}\Bigg)^2\\
  &\geq- C\sigma_n^2\(\frac{\mu_4}{\mu_2}\)^2
  \mathbbm{1}_{\big\{ a_{11}^2+\cdots + a_{nn}^2 >0 \big\}}
  \max_{1\leq i\leq n}\sum_{\substack{1\leq j\leq n}}a_{ij}^2 
\end{align*}
for some $C>0$, since $S_3$ vanishes if $a_{11}=\cdots =a_{nn}=0$.
Applying this and \eqref{eq:S_2} to \eqref{eq:bound2}, we obtain the first  inequality from the assertion. To prove the other one, we use \eqref{eq:bound1} and write 
\begin{align*} 
 d_K\(\frac{Q_n}{\sigma_n}\)\leq \frac{C}{\sigma_n^2}
  \(\mu_4^2\sum_{{1\leq i,j\leq n}}\Bigg(\sum_{\substack{1\leq k\leq n\\k\neq i,j}}a_{ik}a_{kj}\Bigg)^2+\mu_8\sum_{i=1}^na_{ii}^4+\mu_8\sum_{i=1}^n\Bigg(\sum_{\substack{1\leq j\leq n\\i\neq j}}a_{jj}a_{ij}\Bigg)^2\)^{1/2}.
\end{align*}
Next, we bound
\begin{align}\nonumber
&\sum_{{1\leq i,j\leq n}}\Bigg(\sum_{\substack{1\leq k\leq n\\k\neq i,j}}a_{ik}a_{kj}\Bigg)^2=\sum_{{1\leq i,j\leq n}}\Bigg(\sum_{{1\leq k\leq n}}a_{ik}a_{kj}-a_{ii}a_{ij}-a_{ij}a_{jj}\Bigg)^2\\\label{eq:<Tr1}
&\leq \sum_{{1\leq i,j\leq n}}\[2\Bigg(\sum_{{1\leq k\leq n}}a_{ik}a_{kj}\Bigg)^2+4a_{ii}^2a_{ij}^2\]\leq 2{\rm Tr}(A_n^4) +4\sum_{{1\leq i\leq n}}\Bigg(\sum_{{1\leq k\leq n}}a_{ik}^2\Bigg)^2\leq 6{\rm Tr}(A_n^4),
\end{align}
and, by the inequality $ab\leq (a^2+b^2)/2$, 
\begin{align}\nonumber
&\sum_{i=1}^na_{ii}^4+\sum_{i=1}^n\Bigg(\sum_{\substack{1\leq j\leq n\\i\neq j}}a_{jj}a_{ij}\Bigg)^2=\sum_{i=1}^na_{ii}^4+\sum_{i=1}^n\sum_{\substack{1\leq j,k\leq n\\i\neq j}}(a_{ij}a_{kk})(a_{ik}a_{jj})\\\label{eq:<Tr2}
&\leq \sum_{i=1}^na_{ii}^4+\sum_{\substack{1\leq i,j,k\leq n\\i\neq j}}(a_{ij}a_{kk})^2\leq2\sum_{{1\leq i\leq n}}\Bigg(\sum_{{1\leq k\leq n}}a_{ik}^2\Bigg)^2\leq 2{\rm Tr}(A_n^4).
\end{align}
This ends the proof.
\end{Proof}
\noindent
 Contrary to what is stated on page~1590 of \cite{chatterjee}, 
 the conditions
 $\sigma_n^{-2} \sqrt{{\rm Tr}(A_n^4)} \to 0$
 and $\E\[(Q_n/\sigma_n )^4\] \to 3$  
 are not equivalent as $n$ tends to infinity,
  and therefore fourth moment convergence
 is not  sufficient for the central limit theorem to hold
 for quadratic functionals. The next proposition
 clarifies this fact via inequalities between the quantities appearing
 in Theorem \ref{qf}. In the sequel, we let $a\wedge b := \min (a,b)$, $a,b\in \real$. 
 \begin{prop}\label{lem:Tr<Q4+max<Tr}
   There exist absolute constants $C_1, C_2,C_3 >0$ such that 
\begin{align*}
C_1 \frac{\mu_2^4\wedge \tilde\mu_8}{\sigma_n^4} \,{\rm Tr}(A_n^4)&\leq \big|\E\[(Q_n/\sigma_n )^4\]-3\big|+ \frac{\alpha_n^2}{\sigma_n^2} \max_{1\leq i \leq n}\sum_{1\leq j\leq n}a_{ij}^2\\
&\leq  C_2\(\frac{\beta_n^2}{\sigma_n^4}\,{\rm Tr}(A_n^4)+\frac{\alpha_n^2}{\sigma_n^2} \max_{1\leq i \leq n}\sum_{1\leq j\leq n}a_{ij}^2\)\leq C_3\frac{\beta_n^2}{\mu_2^2\sigma_n^2}\,
 \sqrt{ {\rm Tr}(A_n^4) },
\end{align*}
where $\alpha_n, \beta_n$ are as in Theorem \ref{qf}
\end{prop}
\begin{Proof}
  The proof of Theorem~\ref{qf} shows that the
  right hand side of \eqref{r1} is larger than the
  right hand side of \eqref{eq:bound1} up to an absolute multiplicative constant,
  hence we have
$$\big|\E\[(Q_n/\sigma_n )^4\]-3\big|+ \alpha^2_n\sigma_n^2 \max_{1\leq i \leq n}\sum_{1\leq j\leq n}a_{ij}^2\geq C 
\frac{\mu_2^4\wedge \tilde\mu_8}{\sigma_n^4}\Bigg(
  \sum_{i=1}^na_{ii}^4+\sum_{{1\leq i,j\leq n}}\Bigg(\sum_{\substack{1\leq k\leq n\\k\neq i,j}}a_{ik}a_{kj}\Bigg)^2\,\Bigg).$$
  Employing  the inequalities $(a+b)^2\leq 2a^2+2b^2$ and 
   $ab\leq (a^2+b^2)/2$, we get
\begin{align*}
\sum_{{1\leq i,j\leq n}}\Bigg(\sum_{{1\leq k\leq n}}a_{ik}a_{kj}\Bigg)^2&=\sum_{{1\leq i,j\leq n}}\Bigg(\sum_{\substack{1\leq k\leq n\\k\neq i,j}}a_{ik}a_{kj}+a_{ii}a_{ij}+a_{ij}a_{jj}\Bigg)^2\\
&\leq  \sum_{{1\leq i,j\leq n}}\[\Bigg(\sum_{\substack{1\leq i,j\leq n\\k\neq i,j}}a_{ik}a_{kj}\Bigg)^2+8a_{ii}^2a_{ij}^2\]\\
&= \sum_{{1\leq i,j\leq n}}\Bigg(\sum_{\substack{1\leq i,j\leq n\\k\neq i,j}}a_{ik}a_{kj}\Bigg)^2+8\sum_{{1\leq i\leq n}}a_{ii}^2\sum_{\substack{1\leq j\leq n\\j\neq i}}a_{ij}^2+8\sum_{{1\leq i\leq n}}a_{ii}^4\\
&\leq  5\sum_{{1\leq i,j\leq n}}\Bigg(\sum_{\substack{1\leq i,j\leq n\\k\neq i,j}}a_{ik}a_{kj}\Bigg)^2+12\sum_{{1\leq i\leq n}}a_{ii}^4,
\end{align*}  
which gives  the first inequality in the assertion. In order to justify the latter one, we will show
\begin{align}\label{eq:eQ^4<}
\big|\E\[(Q_n/\sigma_n )^4\]-3\big|\leq C\(\frac{\beta_n^2}{\sigma_n^4}\,{\rm Tr}(A_n^4)+\frac{\alpha_n^2}{\sigma_n^2} \max_{1\leq i \leq n}\sum_{1\leq j\leq n}a_{ij}^2\),
\end{align}
  for some $C>0$. Following notation from the proof of Theorem \ref{qf}, we have $\big|\E\[(Q_n/\sigma_n )^4\]-3\big|\leq \(|S_1|+3|S_2-\sigma_n^4|+4|S_3|\)/\sigma_n^4$. By \eqref{eq:S1a}, \eqref{eq:S1b}, \eqref{eq:S_2} and bounding terms from the first  three sums in $S_3$ by $a_{ii}^2a_{ij}^2+a_{ij}^4$ and the last two sums from $S_3$ by
  $$\sum_{i=1}^n\Bigg[\Bigg(\sum_{\substack{1\leq k\leq n\\k\neq i}}a_{ik}^2\Bigg)^2+\Bigg(\sum_{\substack{1\leq j\leq n\\j\neq i}}a_{jj}a_{ij}\Bigg)^2\Bigg],$$ 
  and 
  $$\sum_{\substack{1\leq i,j\leq n\\i\neq j}}\Bigg[a_{ii}^2a^2_{ij}+\Bigg(\sum_{\substack{1\leq k\leq n\\k\neq i,j}}a_{ik}a_{kj}\Bigg)^2\Bigg],$$
  respectively, we arrive at 
  \begin{align*}
  &\big|\E\[(Q_n/\sigma_n )^4\]-3\big|\\
  &\leq 
 C\frac{\beta_n^2}{\sigma_n^4}\Bigg[\sum_{i=1}^na_{ii}^4  +\sum_{{1\leq i,j\leq n}}\Bigg(\sum_{\substack{1\leq k\leq n\\k\neq i,j}}a_{ik}a_{kj}\Bigg)^2+\sum_{i=1}^n\Bigg(\sum_{\substack{1\leq j\leq n\\i\neq j}}a_{jj}a_{ij}\Bigg)^2\Bigg] + \alpha^2_n\sigma_n^2 \max_{1\leq i \leq n}\sum_{1\leq j\leq n}a_{ij}^2,
  \end{align*}
  and \eqref{eq:eQ^4<} follows from \eqref{eq:<Tr1} and \eqref{eq:<Tr2}. Finally, the last bound in the assertion is a consequence of 
  \begin{align*}
 \max_{1\leq i \leq n} \sum_{j=1}^n a_{ij}^2 \leq
 \sqrt{\sum_{i=1}^n\Bigg(\sum_{j=1}^n a_{ij}^2\Bigg)^2}
 \leq \sqrt{{\rm Tr}(A_n^4)},
 \end{align*}
 and 
  $${\rm Tr}(A_n^4)\leq\sum_{i, j=1}^n\Bigg(\sum_{k=1}^n a_{ik}^2\Bigg)\Bigg(\sum_{k=1}^n a_{kj}^2\Bigg)\leq \frac{\sigma_n^4}{\mu_2^4}.$$
\end{Proof}
\noindent 
 Theorem \ref{qf} and Lemma \ref{lem:Tr<Q4+max<Tr} immediately imply
\begin{corollary}
  Assume   $(X_i)_{i\in \N}$ is a fixed i.i.d. sequence with zero means
  and finite $8th$ moments. The following two conditions are equivalent: 
\begin{itemize}
\item[a)] $\E\[(Q_n / \sigma_n )^4\]\longrightarrow 3\ \text{ and }\sigma_n^{-2}\max_{1\leq i \leq n}\sum_{j=1}^na_{ij}^2\longrightarrow0$,
\item[b)]$\sigma_n^{-4}{\rm Tr}(A_n^4)\longrightarrow0,$
\end{itemize}
and they imply $Q_n / \sigma_n \stackrel{\mathcal L}{\longrightarrow} \mathcal N$ with the Kolmogorov rates \eqref{r1} and \eqref{r2}.
\end{corollary}
This extends \eqref{eq:dJ} for any matrix $A_n$ and completes it with the equivalent condition in terms of the trace of $A_n$.

\section*{Acknowledgements}
{G. Serafin was supported by the National Science Centre, Poland, grant no. 2015/18/E/ST1/00239.}

\footnotesize 

\def\cprime{$'$} \def\polhk#1{\setbox0=\hbox{#1}{\ooalign{\hidewidth
  \lower1.5ex\hbox{`}\hidewidth\crcr\unhbox0}}}
  \def\polhk#1{\setbox0=\hbox{#1}{\ooalign{\hidewidth
  \lower1.5ex\hbox{`}\hidewidth\crcr\unhbox0}}} \def\cprime{$'$}

\end{document}